%% file: tifiss_paper_rev.tex
\documentclass[12pt]{article}

\usepackage[a4paper,top=2.5cm,bottom=2.5cm,left=2.5cm,right=2.5cm]{geometry}

\usepackage{amsfonts}
\usepackage{amsmath}
\usepackage{amssymb}
\usepackage{color}
\usepackage{float}
\usepackage{booktabs}			
\usepackage{comment}			
\usepackage{stmaryrd}			
\usepackage{graphicx}			
\usepackage{subfig}				
\usepackage{caption}			
\usepackage{setspace}			
\usepackage{url}
\usepackage{hyperref}			
\usepackage{tikz}
	\usetikzlibrary{calc}
	\usetikzlibrary{positioning}
	\usetikzlibrary{shapes,decorations} 
\usepackage{pgfplots}
	\usepgfplotslibrary{colorbrewer}
	\usepgfplotslibrary{colormaps}
	\pgfplotsset{
		compat = 1.13,												
		every axis/.append style={
		axis x line=box, axis y line=box,							
		tick label style={font={\fontsize{8pt}{12pt}\selectfont}},  
		font={\fontsize{10pt}{12pt}\selectfont}, 					
		xlabel style={font=\footnotesize,yshift=+0.5ex},			
		}
	}
\usepackage[mathscr]{eucal}						

\hypersetup{colorlinks, citecolor=black, filecolor=black, linkcolor=black, urlcolor=black}

\newcommand{\PP}{\mathbb{P}}
\newcommand{\EE}{\mathbb{E}}
\newcommand{\RR}{{\mathbb{R}}}
\newcommand{\NN}{{\mathbb{N}}}
\newcommand{\M}{{\mathcal{M}}}		
\newcommand{\mesh}{\mathcal{T}}		
\newcommand{\xx}{\textbf{x}}

\newcommand{\Span}{\text{span}}                             
\newcommand{\Cov}{\text{Cov}}

\newcommand{\gotI}{\mathscr{I}}		
\newcommand{\gotP}{\mathscr{P}}		
\newcommand{\gotQ}{\mathscr{Q}_\gotP}	
\newcommand{\gotM}{\mathscr{M}}		

\definecolor{myBlue}{RGB}{0,114,200}
\definecolor{lightGreen}{rgb}{0.6,1.0,0.8}	
\definecolor{darkGreen}{RGB}{0,100,0}
\definecolor{myGreen}{RGB}{65,197,116}
\definecolor{myGreenDark}{RGB}{40,160,81}
\definecolor{myGreenDarkTwo}{RGB}{43,120,76}
\definecolor{green3}{HTML}{009688} 
\definecolor{myWater}{RGB}{65,197,216}
\definecolor{Red}{rgb}{1.0,0.0,0.0}
\definecolor{myRed}{rgb}{0.7,0.0,0.0}
\definecolor{myBrown}{rgb}{0.6 0.4 0.2}
\definecolor{myBrownDark}{RGB}{168 114 60}
\definecolor{myOrange}{rgb}{1.0 0.6 0.2}
\definecolor{myLightGray}{RGB}{235,235,235}
\definecolor{myLightGray2}{HTML}{FAFAFA}
\definecolor{myHeavyGray}{HTML}{212121}
\definecolor{myGrayNote}{HTML}{795548}
\definecolor{myCrimson}{RGB}{220,20,60}
\definecolor{myMaroon}{rgb}{0.6,0.2,0.2}
\definecolor{myYellow}{RGB}{255,208,3}
\definecolor{myViolet}{RGB}{153,50,204}
\definecolor{myGray}{RGB}{191,191,191}

\newcommand{\be}{\begin{equation}}
\newcommand{\ee}{\end{equation}}
\newcommand{\bea}{\begin{eqnarray}}
\newcommand{\eea}{\end{eqnarray}}
\newcommand{\beas}{\begin{eqnarray*}}
\newcommand{\eeas}{\end{eqnarray*}}
\newcommand{\ba}{\begin{array}}
\newcommand{\ea}{\end{array}}

\newcommand{\supp}{\text{supp}} 

\newcommand{\G}{\Gamma}

\newcommand{\boldb}{{\mathbf b}}

\newcommand{\boldr}{{\mathbf r}}
\newcommand{\bu}{{\mathbf u}}

\newcommand{\bx}{{\mathbf x}}
\newcommand{\boldy}{{\mathbf y}}

\newcommand{\sfP}{\mathsf{P}}

\newcommand{\enorm}[3][]{#1\vert\hspace*{-.4mm}#1\Vert \, #2 \, #1\Vert\hspace*{-.4mm}#1\vert_{#3}}
\newcommand{\bnorm}[2][]{\enorm[#1]{#2}{}}

\newcommand{\dx}{\mathrm{d}\bx}

\definecolor{otherblue}{rgb}{0,0.3,0.6} 
     \definecolor{maroon}{rgb}{0.6, 0.2, 0.2}
\def\rbl#1{\textcolor{black}{#1}}
\def\abrev#1{\textcolor{black}{#1}}
\def\rev#1{\textcolor{black}{#1}}

\title{\rbl{T-IFISS: a toolbox for adaptive FEM computation}
\thanks{This work was supported by the EPSRC under grants EP/P013317/1 and EP/P013791/1,
\abrev{and was partially supported by The Alan Turing Institute under the EPSRC grant EP/N510129/1}.}}

\author{Alex Bespalov\footnotemark[2] \and Leonardo Rocchi\footnotemark[2] \and David Silvester\footnotemark[3]}

\begin{document}

\date{}
\maketitle

\long\def\symbolfootnote[#1]#2{\begingroup
\def\thefootnote{\fnsymbol{footnote}}\footnote[#1]{#2}\endgroup}

\renewcommand{\thefootnote}{\fnsymbol{footnote}}
\footnotetext[2]{School of Mathematics, University of Birmingham,
Edgbaston, Birmingham B15 2TT, UK ({\tt a.bespalov@bham.ac.uk}, \abrev{\tt leonardo.rocchi@yahoo.it}).}

\footnotetext[3]{\rbl{Department} of Mathematics, University of Manchester,
Oxford Road, Manchester M13 9PL, UK ({\tt d.silvester@manchester.ac.uk}).}

\renewcommand{\thefootnote}{\arabic{footnote}}

\begin{abstract}
\rbl{
T-IFISS is a finite element software package
 for studying finite element solution algorithms for deterministic and parametric elliptic partial
 differential equations. The emphasis is on self-adaptive algorithms with rigorous error control using 
a variety of a posteriori error estimation techniques. The open-source MATLAB 
framework provides a computational laboratory for experimentation and exploration, enabling
users to quickly develop new discretizations and test alternative algorithms.  The package
is also valuable  as  a teaching tool for students who want to learn about
 state-of-the-art finite element methodology.
 }

\end{abstract}

\medskip


\medskip




\section{\rev{Introduction and software overview}}\label{sec:intro}
\setcounter{equation}{0}

Progress in computational mathematics is frequently motivated \abrev{and supported} by the results of 
numerical experiments. The well-established IFISS\footnote{\rev{IFISS is an acronym for
``Incompressible Flow \& Iterative Solver Software''.}} software package~\rev{\cite{ers07, ers14}}
is associated with the monograph~\cite{esw14} and is structured as a stand-alone 
package for studying discretization algorithms for partial differential equations (PDEs) arising in incompressible fluid dynamics.
IFISS is also an established starting point for developing 
code for specialized research applications.\footnote{See the swMATH resource page 
\url{http://swmath.org/software/4398}.}
The package is currently used in universities around the world to enhance the teaching of advanced courses in 
mathematics, computational science and engineering. 
Investigative numerical experiments enable students to develop deduction and 
interpretation skills and are especially useful in helping students  to remember critical ideas in the long term. 

\rev{The \emph{T-IFISS (Triangular IFISS)} toolbox extends the IFISS philosophy and design features
to finite element approximations on \emph{triangular} grids for deterministic and stochastic/parametric
elliptic partial differential equations.}
The emphasis of T-IFISS is on self-adaptive algorithms with rigorous error 
control using a variety of a posteriori error estimation techniques.
\rev{In the same way as for its precursor,
the development of T-IFISS has been motivated by a desire to create a computational laboratory for experimentation
and a tool for reproducible research in computational science;
indeed, the toolbox enables users not only to quickly explore new discretization strategies and test alternative algorithms,
but to replicate, validate and independently verify computational results.}
\rev{Thus, in this article, instead of giving a comprehensive technical description of the software,
we aim to highlight and document
those features of T-IFISS that are not available in its precursor package IFISS;
we will illustrate these features with a series of case studies that demonstrate the efficiency of the software
and its utility as a research and teaching~tool.
}

All the test problems that are built into the current version of T-IFISS\footnote{\rbl{T-IFISS version~1.2
was released in February 2019 and runs under MATLAB or Octave. It can be downloaded 
from \url{http://www.manchester.ac.uk/ifiss/tifiss.html} and is 
compatible with Windows, Linux and MacOS computers.}}
are associated with steady-state diffusion equations
with anisotropic (or uncertain) conductivity coefficients together with \abrev{mixed (Dirichlet and Neumann)} boundary conditions.
These PDE problems may be solved on general polygonal domains, including slit domains and domains with holes.
It is worth pointing out that many of these (deterministic) test problems 
could also be solved using general purpose finite element software packages like
deal.II~\cite{BangerthHartmannKanschat2007}, 
DUNE~\cite{dune2016}, 
FEniCS~\cite{LoggMardalEtAl2012a}, 
and FreeFEM~\cite{Hecht_12_NDF}. 
The attraction of the T-IFISS environment is the ease with which one can test 
the alternative error estimation, marking, and refinement strategies, develop new strategies, and
extend functionality to new problem classes as well as new approximations and solution algorithms.
\rev{The modular structure of the code allows one
to examine interactions between different components in order to achieve
favourable convergence properties or minimize the associated computational cost.
Finally, the use of the high-level MATLAB syntax ensures readability, accessibility, 
and quick visualization of the solution, error estimators, and finite element meshes.}
These features \rev{are invaluable in teaching but} are not so explicit in other adaptive finite  element method (FEM) packages like 
ALBERTA~\cite{SchmidtS_05_DAF}, 
PLTMG~\cite{Bank_98_PLTMG}, 
and p1afem~\cite{FunkenPW_11_EIA}. 
\rev{A unique feature of T-IFISS is that it can be used to solve
parameter-dependent elliptic PDE problems stemming from uncertainty quantification models.
This facility, called {\it stochastic} T-IFISS~\cite{BespalovR_stoch_tifiss}, develops the idea of
adaptive stochastic Galerkin FEM and provides an effective alternative to traditional sampling
methods commonly used for such problems. The only software  packages that we know of
that have a similar  capability are ALEA~\cite{alea} and  SGLib~\cite{sglib}.}

\rev{
The following is a brief overview of the functionality available in T-IFISS (including Stochastic T-IFISS)
with links to its directory structure
(we refer to~\cite[Appendix~B]{Rocchi_thesis} for more details including the description of the associated data structures):
%
\begin{itemize}
\item
computing Galerkin solutions with $P_1$ and $P_2$ approximations over the specified triangular mesh
(directory {\tt /diffusion});
\item
computing stochastic Galerkin FEM solutions with spatial $P_1$ and $P_2$ approximations
over the specified triangular mesh and for the specified polynomial approximation on the parameter domain;
the functionality here includes effective iterative solution of very large linear systems stemming from stochastic Galerkin FEM
(directory {\tt /stoch\_diffusion});
\item
estimating the energy error in the computed Galerkin solution and
the error in a quantity of interest derived from this solution
(directories {\tt /diffusion\_error} and {\tt /goafem}, respectively);
\item
estimating the energy error in the computed stochastic Galerkin solution and
the error in a quantity of interest derived from this solution
(directories {\tt /stoch\_diffusion} and {\tt /stoch\_diffusion/stoch\_goafem}, respectively);
\item
adaptive refinement of Galerkin approximations, including local adaptive mesh refinement
and, in the case of stochastic/parametric problems, adaptive enrichment of multivariable
polynomial approximations
(directories {\tt /diffusion\_adaptive} and
{\tt /stoch\_diffusion/stoch\_diffusion\_adapt});
\item
visualization of Galerkin solutions, goal functionals, error estimators, finite element meshes;
plotting convergence history (the error estimates against the number of degrees of freedom) 
(directory {\tt /graphs});
\end{itemize}
}

\rev{The rest of the article is organized as follows.
In the next section, we recall the main ingredients of adaptive FEM and describe their implementation in T-IFISS;
we illustrate this  with two case studies (solving the diffusion equation with strongly anisotropic coefficient
and computing a harmonic function in the L-shaped domain).
Section~\ref{sec:goafem} is focused on implementing the goal-oriented error estimation and adaptivity;
the case study here demonstrates the capability of the software to approximate the value of a singular solution
to the diffusion problem at a fixed point away from the singularity.
In section~\ref{sec:asgfem}, we discuss the implementation of stochastic Galerkin FEM including the
associated error estimation and adaptivity.
The efficiency of our adaptive algorithm is demonstrated with a representative case study
(solving the diffusion equation with parametric coefficient over the L-shaped domain).
Some extensions and future developments of the toolbox are briefly outlined in Section~\ref{sec:extensions}.
}

\section{\rbl{Adaptive finite  element methods}}\label{sec:adaptivity}
\setcounter{equation}{0}

Adaptive finite element approximations to solutions of partial differential equations
are typically computed by iterating the following loop of four component modules:
\be \label{eq:modules}
    \mathsf{SOLVE}\ \Longrightarrow\ \mathsf{ESTIMATE}\ 
\Longrightarrow\ \mathsf{MARK}\ \Longrightarrow\ \mathsf{REFINE}.
\ee

\rev{
In this section we first describe the implementation of the four components in~\eqref{eq:modules} within T-IFISS.
This is followed by two case studies that highlight some of the features 
of the  toolbox and demonstrate its utility and its efficiency.
}


\subsection{\rev{Main ingredients of adaptive FEM and their implementation in T-IFISS}} \label{subsec:adaptivity:ingredients}

\rev{Module} $\mathsf{SOLVE}$.
For a given PDE problem, the Galerkin solution \rbl{defined on a specific}
\rev{(structured or unstructured\footnote{\rev{In the MATLAB release of T-IFISS,
unstructured meshes are generated using the DistMesh package~\cite{PerssonS_04_SMG}.}})} triangular mesh is computed
by solving the linear system \rbl{associated with the} Galerkin projection of the variational 
formulation of the problem onto the corresponding finite element space.
Two types of ($C^0$) finite element approximations are implemented:  piecewise linear ($P_1$) and 
piecewise quadratic ($P_2$). For either of these approximations, fast computation of the entries in the
 stiffness matrix and the load vector is \rbl{achieved by vectorizing the calculations over all the  elements.
 When  solving a deterministic problem the resulting linear equation system is solved using the 
 highly optimized sparse direct solver (UMFPACK) that is built  into MATLAB and Octave}.\footnote{\rbl{We would 
 almost certainly use an iterative solver  preconditioned with an algebraic multigrid V-cycle  
 if we were trying to solve the same  PDE problem in three spatial dimensions.}}

\smallskip

\rev{Module} $\mathsf{ESTIMATE}$.
\rbl{The purpose of this module is two-fold.}
First, it computes local error indicators that provide information about the distribution
of estimated local errors in the computed Galerkin solution;
the error indicators may be associated with elements or edges of the underlying triangulation.
Second, it computes an estimate of the (total) energy error in the Galerkin solution.
\rbl{This estimate} is used to decide whether the stopping tolerance is met.
T-IFISS offers a choice of \rbl{the following} three error estimation strategies (EES) for
 {\it linear} ($P_1$) approximation. 

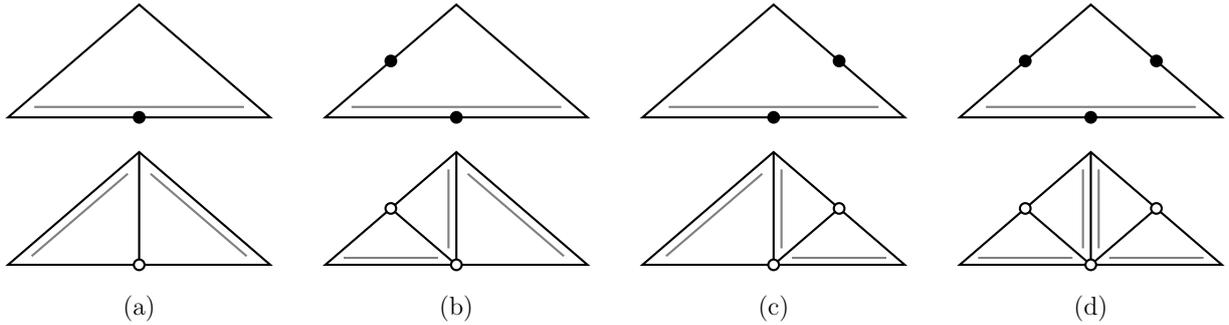
\begin{figure}[!t]
\centering
\begin{tikzpicture}[thick,scale=1.15]
%
\node[scale=0.81] at (1.5,-0.5) {(a)};
%
\draw (0,0) -- (3,0) -- (1.5,1.3) -- cycle;
\draw (0,1.7) -- (3,1.7) -- (1.5,3) -- cycle;
\draw (1.5,0) -- (1.5,1.3);
\draw[color=gray] (0.3,1.82) -- (2.7,1.82);
\draw[color=gray] (0.27,0.1) -- (1.37,1.05);
\draw[color=gray] (2.73,0.1) -- (1.63,1.05);
\node[circle,draw=black, fill=black, inner sep=0pt,minimum size=4pt, scale=1] (b) at (1.5,1.7) {};
\node[circle,draw=black, fill=white, inner sep=0pt,minimum size=4pt, scale=1] (b) at (1.5,0) {};
\end{tikzpicture}
\hfill
\begin{tikzpicture}[thick,scale=1.15]
%
\node[scale=0.81] at (1.5,-0.5) {(b)};
%
\draw (0,0) -- (3,0) -- (1.5,1.3) -- cycle;
\draw (0,1.7) -- (3,1.7) -- (1.5,3) -- cycle;
\draw (1.5,0) -- (1.5,1.3);
\draw (1.5,0) -- (0.75,0.65);
\draw[color=gray] (0.3,1.82) -- (2.7,1.82);		
\draw[color=gray] (2.73,0.1) -- (1.63,1.05);
\draw[color=gray] (0.21,0.08) -- (1.29,0.08);	
\draw[color=gray] (1.41,0.19) -- (1.41,1.1);	
\node[circle,draw=black, fill=black, inner sep=0pt,minimum size=4pt, scale=1] (b) at (1.5,1.7) {};
\node[circle,draw=black, fill=black, inner sep=0pt,minimum size=4pt, scale=1] (b) at (0.75,2.35) {};
\node[circle,draw=black, fill=white, inner sep=0pt,minimum size=4pt, scale=1] (b) at (1.5,0) {};
\node[circle,draw=black, fill=white, inner sep=0pt,minimum size=4pt, scale=1] (b) at (0.75,0.65) {};
\end{tikzpicture}
\hfill
\begin{tikzpicture}[thick,scale=1.15]
%
\node[scale=0.81] at (1.5,-0.5) {(c)};
%
\draw (0,0) -- (3,0) -- (1.5,1.3) -- cycle;
\draw (0,1.7) -- (3,1.7) -- (1.5,3) -- cycle;
\draw (1.5,0) -- (1.5,1.3);
\draw (1.5,0) -- (2.25,0.65);
\draw[color=gray] (0.3,1.82) -- (2.7,1.82);   
\draw[color=gray] (0.27,0.1) -- (1.37,1.05);  
\draw[color=gray] (1.71,0.08) -- (2.79,0.08); 
\draw[color=gray] (1.59,0.19) -- (1.59,1.1);  
\node[circle,draw=black, fill=black, inner sep=0pt,minimum size=4pt, scale=1] (b) at (1.5,1.7) {};
\node[circle,draw=black, fill=black, inner sep=0pt,minimum size=4pt, scale=1] (b) at (2.25,2.35) {};
\node[circle,draw=black, fill=white, inner sep=0pt,minimum size=4pt, scale=1] (b) at (1.5,0) {};
\node[circle,draw=black, fill=white, inner sep=0pt,minimum size=4pt, scale=1] (b) at (2.25,0.65) {};
\end{tikzpicture}
\hfill
\begin{tikzpicture}[thick,scale=1.15]
%
\node[scale=0.81] at (1.5,-0.5) {(d)};
%
\draw (0,0) -- (3,0) -- (1.5,1.3) -- cycle;
\draw (0,1.7) -- (3,1.7) -- (1.5,3) -- cycle;
\draw (1.5,0) -- (1.5,1.3);
\draw (1.5,0) -- (2.25,0.65);
\draw (1.5,0) -- (0.75,0.65);
\draw[color=gray] (0.3,1.82) -- (2.7,1.82);   
\draw[color=gray] (1.71,0.08) -- (2.79,0.08); 
\draw[color=gray] (1.59,0.19) -- (1.59,1.1);  
\draw[color=gray] (0.21,0.08) -- (1.29,0.08);  
\draw[color=gray] (1.41,0.19) -- (1.41,1.1);  
\node[circle,draw=black, fill=black, inner sep=0pt,minimum size=4pt, scale=1] (b) at (1.5,1.7) {};
\node[circle,draw=black, fill=black, inner sep=0pt,minimum size=4pt, scale=1] (b) at (2.25,2.35) {};
\node[circle,draw=black, fill=white, inner sep=0pt,minimum size=4pt, scale=1] (b) at (1.5,0) {};
\node[circle,draw=black, fill=white, inner sep=0pt,minimum size=4pt, scale=1] (b) at (2.25,0.65) {};
\node[circle,draw=black, fill=black, inner sep=0pt,minimum size=4pt, scale=1] (b) at (0.75,2.35) {};
\node[circle,draw=black, fill=white, inner sep=0pt,minimum size=4pt, scale=1] (b) at (0.75,0.65) {};
\end{tikzpicture}
\caption{ 
NVB bisections:
(a) one, (b)-(c) two, and (d) three bisections of the edges of the triangular element.
Double lines indicate reference edges, black dots indicate the edges to be bisected,
and white dots indicate the newest vertices.
}
\label{fig:nvb:ref}
\end{figure}

(EES1) \rbl{is a local hierarchical error estimator computed via} a standard element residual 
technique (see~\cite[Section~3.3]{AinsworthOden00})
using either piecewise linear or piecewise quadratic bubble functions
over subelements obtained by uniform refinements\footnote{The default uniform refinement in T-IFISS
is by three bisections (see Figure~\ref{fig:nvb:ref}(d)). However, there is an option to switch to the 
so-called \emph{red} uniform refinement
(i.e., the one obtained by connecting the edge midpoints of each triangle);
this can be done by setting {\tt subdivPar = 1} within the function {\tt adiff\_adaptive\_main.m}.}.
The local error indicators in this case are computed elementwise by solving $3 \times 3$ linear systems
(this calculation is vectorized over elements).
The total error estimate is calculated as the $\ell_2$-norm of the vector of local error indicators.

(EES2) \rbl{is} a global hierarchical estimator (see~\cite{BankW_85_SAE},~\cite[Section~5]{AinsworthOden00}) using
piecewise linear bubble functions corresponding to the uniform refinement of the original triangulation.
Note that the implementation of this strategy requires solving a sparse linear system associated with 
a global residual problem.
The localizations of the estimator (to either the elements (default option) or the 
edges\footnote{More precisely, the interior edges and the edges associated with those parts of 
the boundary where the Neumann and non-homogeneous Dirichlet boundary conditions are prescribed.}
of triangulation) gives two types of local error indicators in this EES.

(EES3) \rbl{is}  a two-level error estimate employing piecewise linear bubble functions associated
with edge midpoints (of the original triangulation); see~\cite{MSW98, ms99}.
In this case, it is natural to choose local error indicators associated with edges
(this is the default choice in (EES3)).
However, the choice of elementwise error indicators is also offered as an option;
these are computed for each interior element from three corresponding edge indicators
(or from two indicators for the elements with an edge on the Dirichlet part of the boundary).

\rbl{There is \abrev{currently} no flexibility \abrev{with choosing the EES} when using  {\it quadratic} ($P_2$) approximation:}
the local hierarchical \abrev{error} estimation strategy (EES1) is employed with piecewise quartic bubble functions.
\rbl{More specifically,} the local error indicators are computed elementwise by solving $9 \times 9$ linear systems
(again, the calculation is vectorized over elements), and
the total error estimate is calculated as the $\ell_2$-norm of the vector of local error indicators.

\smallskip

\rev{Module} $\mathsf{MARK}$.
\rbl{In this module} the elements (or edges) with largest error indicators are selected (i.e., marked) for refinement.
Two popular marking strategies \rbl{are currently implemented}:
the {\it maximum} strategy and the {\it D\"orfler} strategy (also referred to as the equilibration or bulk chasing strategy).

Let $\{\beta(s);\; s \in \mathcal{S}\}$ denote the set of error indicators associated with the elements of the set 
$\mathcal{S}$ (e.g., $\mathcal{S}$ can be the set of edges or elements of the triangulation).
In the maximum marking strategy, that dates back to~\cite{bv84},
the element $s \in \mathcal{S}$ is marked if the associated error indicator $\beta(s)$ is larger
than a fixed proportion of the maximum among all error indicators.
Specifically, for a given marking (or, threshold) parameter $\theta \in [0, 1]$,
this strategy returns a minimal subset $\mathcal{M} \subseteq \mathcal{S}$ of marked elements such that
\be \label{marking:max}
      \beta(s) \ge \theta \max_{s \in \mathcal{S}} \beta(s)\qquad \forall\, s \in \mathcal{M}.
\ee
Note that in this strategy, smaller values of $\theta$ lead to larger subsets $\mathcal{M}$.

In the D\"orfler marking strategy, that was originally introduced in~\cite{doerfler},
sufficiently many elements are marked such that the combined contribution
of the corresponding error indicators is larger than a fixed proportion of the total error estimate.
More precisely, given a marking parameter $\theta \in (0, 1]$,
this strategy builds a subset $\mathcal{M} \subseteq \mathcal{S}$ of minimal cardinality such that
$\{\beta(s);\; s \in \mathcal{M}\}$ is the set of $\#\mathcal{M}$ largest error indicators and
\be \label{marking:doerfler}
      \sum\limits_{s \in \mathcal{M}} \beta(s)^2 \ge \theta \sum\limits_{s \in \mathcal{S}} \beta(s)^2.
\ee
Here, smaller values of $\theta$ lead to smaller subsets $\mathcal{M}$.

\smallskip

\rev{Module} $\mathsf{REFINE}$.
Given the set of marked elements (or, marked edges) that is obtained by employing one of the
above marking strategies, local adaptive mesh refinement is performed in T-IFISS by implementing
the longest edge bisection (LEB) strategy---a variant of the newest vertex bisection~(NVB) method
(we refer to~\cite{s72Phd, rivara84, b91, k95, s08, nv12} for theoretical and implementational aspects 
of NVB refinements,
as well as to \cite{m89} for an overview and comparison of NVB  with other mesh-refinement techniques).
In this method, a \emph{reference edge} is designated for each triangle $T$
(for the coarsest mesh, this is always the longest edge of $T$), and $T$ is bisected by halving the 
reference edge; see Figure~\ref{fig:nvb:ref}(a).
This introduces two new elements, the sons of $T$, for which reference edges are selected\footnote{In the NVB 
method, the reference edges of the sons are the edges opposite to the new vertex, whereas in the LEB method,
 the reference edge is always the longest edge of the element. Note, however, that for structured triangulations 
 of square, L-shaped, and crack domains, both methods result in identical refinement patterns.}.
A recursive application of this procedure leads to a conforming mesh,
where one, two, or three bisections of the triangle $T$ may be performed;
see Figure~\ref{fig:nvb:ref}(a)--(d).
The refinement by three bisections (see Figure~\ref{fig:nvb:ref}(d)) is called the \emph{bisec3} refinement.
Refining all elements of the given mesh by three bisections results in a (conforming) uniform \emph{bisec3}
 refinement of this mesh.

It is important to emphasize that NVB iteratively refines individual elements by bisecting some (or all) of their edges.
Therefore, either the set of marked elements or the set of marked edges can be used as an input
to the NVB-based mesh-refinement routine.
Furthermore, NVB refinements lead to nested (Lagrange) finite element spaces (see~\cite[p.179]{nv12})---an important
ingredient in the proof of the contraction property for adaptive finite element approximations, see~\cite[Section~5]{nv12}
(note that nestedness is not guaranteed for other mesh-refinement techniques, such as
red--green or red--green--blue~refinements).

As many other modules in the toolbox, the mesh-refinement routine in T-IFISS exploits MATLAB's vectorization features.
More precisely, once the set of marked elements or edges is returned by the module $\mathsf{MARK}$, 
the mesh-refinement routine identifies the subsets of elements where one, two, or three 
bisections should be performed (see Figure~\ref{fig:nvb:ref}) and then the elements in the 
three separate subsets are refined simultaneously.

\subsection{\rev{Numerical case studies}} \label{subsec:adaptivity:numerics}

\rev{Let us demonstrate the performance of adaptive finite element routines in T-IFISS with two test examples.
When doing this, we will illustrate some of the ingredients of adaptive FEM described in the previous subsection.}

\smallskip

\begin{figure}[b!] 
\centering
\subfloat[][]
{
\includegraphics[scale=1.29]{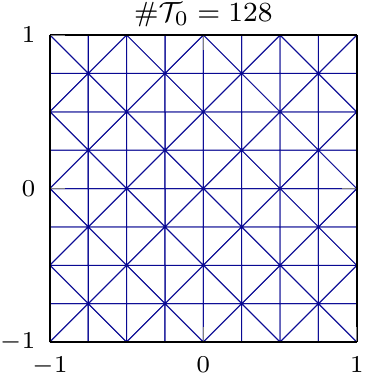}} 
\hspace{2cm}
\subfloat[][] 
{
\includegraphics[scale = 0.94]{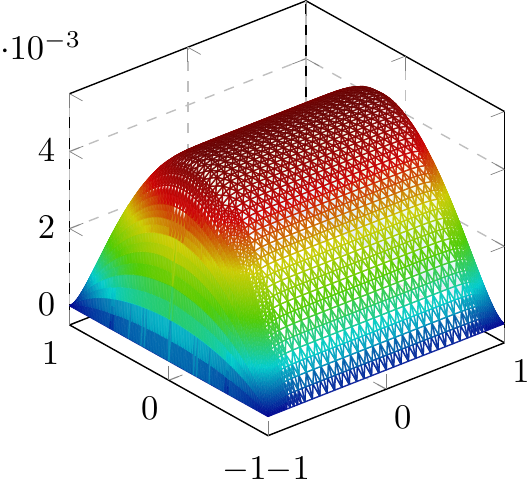}}
\caption{
\rev{Example 1: (a) initial coarse mesh; 
(b) Galerkin solution to problem~\eqref{eq:ex1:problem}.}
}
\label{fig:ex1:sol}
\end{figure}

{\bf Example 1.} \
Let $D = (-1,1)^2$ be the square domain.
We consider the diffusion equation with \rbl{a strongly} anisotropic coefficient, \rbl{together with  a 
constant source function and a homogeneous Dirichlet boundary condition:}
\be \label{eq:ex1:problem}
\begin{aligned} 
           -\nabla \cdot (A \nabla u(\bx) ) & = 1, \qquad \bx = (x_1,x_2) \in D,\\[3pt]
           u(\bx) & = 0, \qquad \bx \in \partial D,
\end{aligned} 
\ee
where $A = \left[\begin{smallmatrix} 1 & 0\\ 0 & 100 \end{smallmatrix} \right]$. \rbl{ We solve this problem 
using  $P_1$ approximation. The solution  is depicted in Figure~\ref{fig:ex1:sol}(b) and} 
exhibits sharp gradients within the boundary layers  along the edges $x_1 = \pm 1$  of the domain.

\begin{figure}[t!]
\centering
\subfloat[][] 
{
\begin{tikzpicture}
\pgfplotstableread{data/exp1/P1-hierarch-theta05-elem.dat}{\first}
\begin{loglogaxis}
[
width = 6.2cm, height=6.0cm,
title={(EES1)}, 
title style={font={\fontsize{9pt}{12pt}\selectfont},yshift=-1.0ex},
xlabel={degrees of freedom}, 
xlabel style={font={\fontsize{8pt}{12pt}\selectfont}},
ylabel={error estimate},
ylabel style={font={\fontsize{8pt}{12pt}\selectfont},yshift=-1.0ex},
ymajorgrids=true, xmajorgrids=true, grid style=dashed, 
xmin = (3)*10^(1), 	xmax = (4)*10^(4),
ymin = (5)*10^(-4),	ymax = (1.5)*10^(-1),
font={\fontsize{8pt}{12pt}\selectfont}, 
legend style={legend pos=south west, legend cell align=left, fill=none, draw=none, 
font={\fontsize{9pt}{12pt}\selectfont}}
]
\addplot[blue,mark=triangle,mark size=3.5pt]	table[x=dofs, y=error]{\first};
\addplot[black,solid,domain=10^(1):10^(7)]    { (0.6)*x^(-0.5) };
\node at (axis cs:2e3,1e-2)	[anchor=south west]{$\mathcal{O}(N^{-0.5})$};
\end{loglogaxis}
\end{tikzpicture}
}
\hfill
\subfloat[][] 
{
\begin{tikzpicture}
\pgfplotstableread{data/exp1/P1-global-hierach-theta05-elem.dat}{\first}
\begin{loglogaxis}
[
width = 6.2cm, height=6.0cm,
title={(EES2)}, 
title style={font={\fontsize{9pt}{12pt}\selectfont},yshift=-1.0ex},
xlabel={degrees of freedom}, 
xlabel style={font={\fontsize{8pt}{12pt}\selectfont}},
ymajorticks=false,
ylabel style={font={\fontsize{8pt}{12pt}\selectfont},yshift=-1.0ex},
ymajorgrids=true, xmajorgrids=true, grid style=dashed, 
xmin = (3)*10^(1), 	xmax = (4)*10^(4),
ymin = (5)*10^(-4),	ymax = (1.5)*10^(-1),
font={\fontsize{8pt}{12pt}\selectfont}, 
legend style={legend pos=south west, legend cell align=left, fill=none, draw=none, 
font={\fontsize{9pt}{12pt}\selectfont}}
]
\addplot[blue,mark=triangle,mark size=3.5pt]		table[x=dofs, y=error]{\first};
\addplot[black,solid,domain=10^(1):10^(7)]    { (0.6)*x^(-0.5) };
\node at (axis cs:2e3,1e-2)	[anchor=south west]{$\mathcal{O}(N^{-0.5})$};
\end{loglogaxis}
\end{tikzpicture}
}
\hfill
\subfloat[][] 
{
\begin{tikzpicture}
\pgfplotstableread{data/exp1/P1-twolevel-theta05-elem.dat}{\first}
\begin{loglogaxis}
[
width = 6.2cm, height=6.0cm,
title={(EES3)}, 
title style={font={\fontsize{9pt}{12pt}\selectfont},yshift=-1.0ex},
xlabel={degrees of freedom}, 
xlabel style={font={\fontsize{8pt}{12pt}\selectfont}},
ymajorticks=false,
ylabel style={font={\fontsize{8pt}{12pt}\selectfont},yshift=-1.0ex},
ymajorgrids=true, xmajorgrids=true, grid style=dashed, 
xmin = (3)*10^(1), 	xmax = (4)*10^(4),
ymin = (5)*10^(-4),	ymax = (1.5)*10^(-1),
font={\fontsize{8pt}{12pt}\selectfont}, 
legend style={legend pos=south west, legend cell align=left, fill=none, draw=none, 
font={\fontsize{9pt}{12pt}\selectfont}}
]
\addplot[blue,mark=triangle,mark size=3.5pt]	table[x=dofs, y=error]{\first};
\addplot[black,solid,domain=10^(1):10^(7)]    { (0.6)*x^(-0.5) };
\node at (axis cs:2e3,1e-2)	[anchor=south west]{$\mathcal{O}(N^{-0.5})$};
\end{loglogaxis}
\end{tikzpicture}
}
\caption{
Example 1: 
error estimates at each iteration of the adaptive algorithm
employing the error estimation strategies (EES1)--(EES3) and using the \emph{element-based} D\"orfler marking 
with $\theta = 0.5$.
Here, $N$ denotes the number of degrees of freedom.
}
\label{fig:ex1:exp1:converge}
\end{figure}
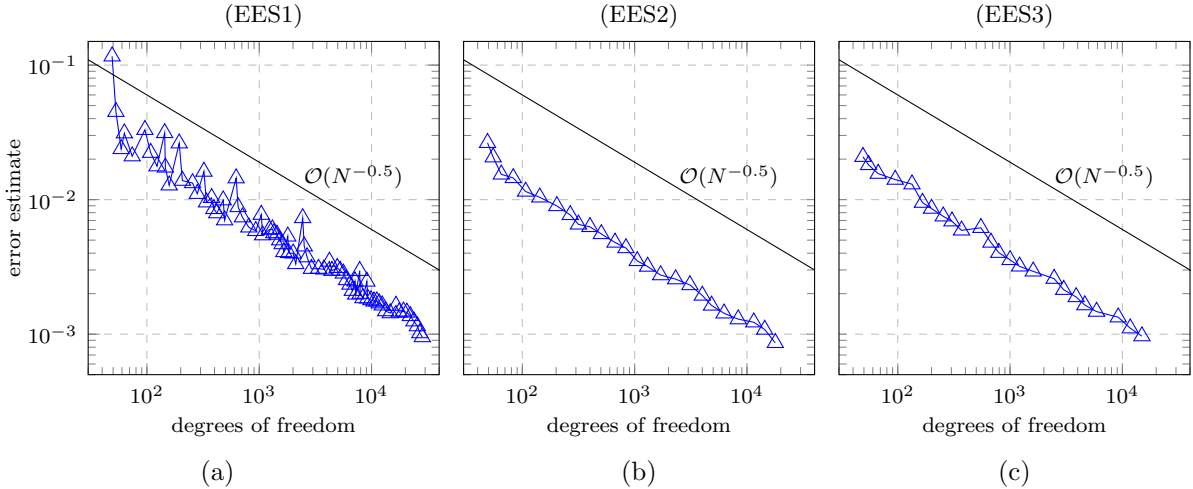

\begin{figure}[b!] 
\centering
\subfloat[][] 
{
\includegraphics[scale = 1.25]{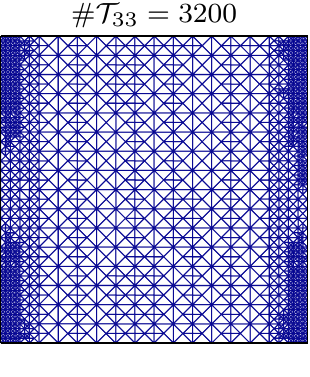}} 
\hfill 
\subfloat[][]
{
\includegraphics[scale = 1.25]{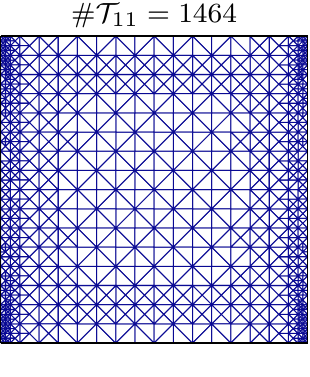}}
\hfill
\subfloat[][] 
{
\includegraphics[scale = 1.25]{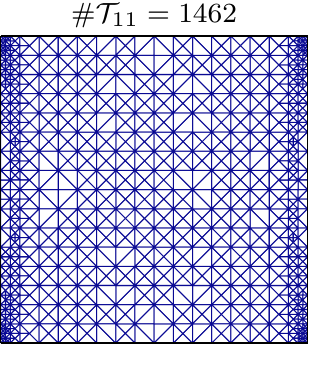}}
\caption{
Example 1: 
\rev{locally refined meshes produced by the adaptive algorithm employing the error estimation strategies 
(EES1)--(EES3) and
using the \emph{element-based} D\"orfler marking with $\theta=0.5$.  \rbl{The header $\#\mathcal{T}_{\ell}$ refers to
 the number of elements in the mesh at step $\ell$ of the adaptive process.}}
}
\label{fig:ex1:exp1:mesh}
\end{figure}

In \rbl{our  first experiment for this problem, we ran} the adaptive FEM algorithm  three times,
employing three different error estimation strategies (EES1)--(EES3) with piecewise linear bubble functions 
as described \rbl{earlier}. In each case, we started with the \rbl{coarse grid of 128 elements  shown in}
Figure~\ref{fig:ex1:sol}(a), \rbl{fixed the stopping tolerance to} {\tt tol = 1e-3}
and employed the {element-based} D\"orfler marking strategy~\eqref{marking:doerfler} with $\theta = 0.5$.
The results of computations are presented in Figures~\ref{fig:ex1:exp1:converge},~\ref{fig:ex1:exp1:mesh}
and in Table~\ref{tbl:ex1:exp1:data}.

Figure~\ref{fig:ex1:exp1:converge} shows convergence plots for the energy error estimates computed 
via each of the three strategies.
In Figure~\ref{fig:ex1:exp1:mesh}, 
we plot the locally refined meshes generated by the adaptive algorithm
in each of the three cases (for some intermediate tolerance).
In Table~\ref{tbl:ex1:exp1:data}, we have collected the data on iteration counts and mesh refinements for each run of the algorithm.
It is evident from these results that \rbl{using (EES1)} leads to unstable reductions in the estimated errors and,
as a consequence, to a large number of iterations and an over-refined final mesh.
\rev{This is because for problem~\eqref{eq:ex1:problem} with constant coefficients,
the elementwise interior residuals for $P_1$ approximations have zero contributions to 
the associated error estimator.}
In contrast, the adaptive algorithms employing (EES2) and (EES3) lead to essentially monotonic decay of the error estimates;
\rbl{the number of iterations needed to reach the stopping tolerance is nearly the same
and similar mesh refinement patterns are generated in both cases.}
We note that for all three strategies, the error estimates decay with an optimal rate of $\mathcal{O}(N^{-0.5})$,
where $N$ is the number of degrees of freedom~(d.o.f.).

\begin{table}[t!]
\setlength\tabcolsep{30pt} 
\begin{center} 
\footnotesize{
\renewcommand{\arraystretch}{1.15}
\begin{tabular}{l   l  l   l  l    l  l   }
\toprule
\multicolumn{7}{c}{Adaptive FEM with \emph{element-based} D\"orfler marking ($\theta=0.5$); \ {\tt tol=1e-3}} \\
\midrule
&\multicolumn{2}{l}{(EES1)}	
&\multicolumn{2}{l}{(EES2)}
&\multicolumn{2}{l}{(EES3)}\\[-1pt]
\midrule 
$L$					&\multicolumn{2}{l}{$77$}
					&\multicolumn{2}{l}{$25$}	
					&\multicolumn{2}{l}{$24$}\\[-3pt]

$\eta_L$			&\multicolumn{2}{l}{$9.520$\text{e-}$04$}	
					&\multicolumn{2}{l}{$8.619$\text{e-}$04$}	
					&\multicolumn{2}{l}{$9.668$\text{e-}$04$}\\[-3pt]
					
$\#\mesh_L$			&\multicolumn{2}{l}{57,347}
					&\multicolumn{2}{l}{36,390}
					&\multicolumn{2}{l}{30,498}\\[-3pt]

$\abrev{n_L}$				&\multicolumn{2}{l}{28,296}	
					&\multicolumn{2}{l}{17,875}	
					&\multicolumn{2}{l}{14,934}\\[-1pt]
\bottomrule 
\end{tabular}
\caption{
\rbl{\abrev{Example 1: output} when solving  \eqref{eq:ex1:problem} using alternative error estimation strategies;}
$L$ denotes the total number of iterations, $\eta_L$, $\#\mathcal{T}_L$, and $\abrev{n_L}$ refer to
the final error estimate, the number of elements in the final mesh, and the number of \rbl{degrees of freedom}
at the \rbl{final}  iteration, respectively.
}
\label{tbl:ex1:exp1:data}       
} 
\end{center}                                                                   
\end{table}

\begin{figure}[h!] 
\centering
\subfloat[][] 
{
\begin{tikzpicture}
\pgfplotstableread{data/exp1/P1-twolevel-theta05-edge.dat}{\first}
\begin{loglogaxis}
[
width = 6.2cm, height=6.0cm,
title={(EES3)}, 
title style={font={\fontsize{9pt}{12pt}\selectfont},yshift=-1.0ex},
xlabel={degrees of freedom}, 
xlabel style={font={\fontsize{8pt}{12pt}\selectfont}},
ylabel={error estimate},
ylabel style={font={\fontsize{8pt}{12pt}\selectfont},yshift=-1.0ex},
ymajorgrids=true, xmajorgrids=true, grid style=dashed, 
xmin = (3)*10^(1), 	xmax = (4)*10^(4),
ymin = (5)*10^(-4),	ymax = (4)*10^(-2),
font={\fontsize{8pt}{12pt}\selectfont}, 
legend style={legend pos=south west, legend cell align=left, fill=none, draw=none, 
font={\fontsize{9pt}{12pt}\selectfont}}
]
\addplot[blue,mark=triangle,mark size=3.5pt]	table[x=dofs, y=error]{\first};
\addplot[black,solid,domain=10^(1):10^(7)]    { (0.3)*x^(-0.5) };
\node at (axis cs:2e3,6e-3)	[anchor=south west]{$\mathcal{O}(N^{-0.5})$};
\end{loglogaxis}
\end{tikzpicture}
} 
\hspace{2cm}
\subfloat[][]
{
\begin{tikzpicture}
\pgfplotstableread{data/exp1/P1-twolevel-theta05-edge.dat}{\first}
\begin{semilogxaxis}
[
width = 6.2cm, height=6.0cm,
title={(EES3)}, 
title style={font={\fontsize{9pt}{12pt}\selectfont},yshift=-1.0ex},
xlabel={degrees of freedom}, 
xlabel style={font={\fontsize{8pt}{12pt}\selectfont}},
ylabel={effectivity index},
ylabel style={font={\fontsize{8pt}{12pt}\selectfont},yshift=-1.0ex},
ymajorgrids=true, xmajorgrids=true, grid style=dashed, 
xmin = (3)*10^(1), 	xmax = (3)*10^4,
ymin = 0.45,		ymax = 0.7, 
ytick={0.5,0.6,0.7,0.8,0.9,1.0},
font={\fontsize{8pt}{12pt}\selectfont}, 
legend style={legend pos=south west, legend cell align=left, fill=none, draw=none, 
font={\fontsize{9pt}{12pt}\selectfont}}
]
\addplot[blue,mark=triangle,mark size=3.5pt]	table[x=dofs, y=effindices]{\first};
\addplot[black,solid,domain=10^(1):10^(7)]    { (0.6)*x^(-0.5) };
\node at (axis cs:2e3,1e-2)	[anchor=south west]{$\mathcal{O}(N^{-0.5})$};
\end{semilogxaxis}
\end{tikzpicture}
}
\caption{
Example 1:
(a) error estimates at each iteration of the adaptive algorithm employing (EES3) and
the \emph{edge-based} D\"orfler marking with $\theta = 0.5$,
(b) the associated effectivity~indices.
}
\label{fig:ex1:exp2}
\end{figure}
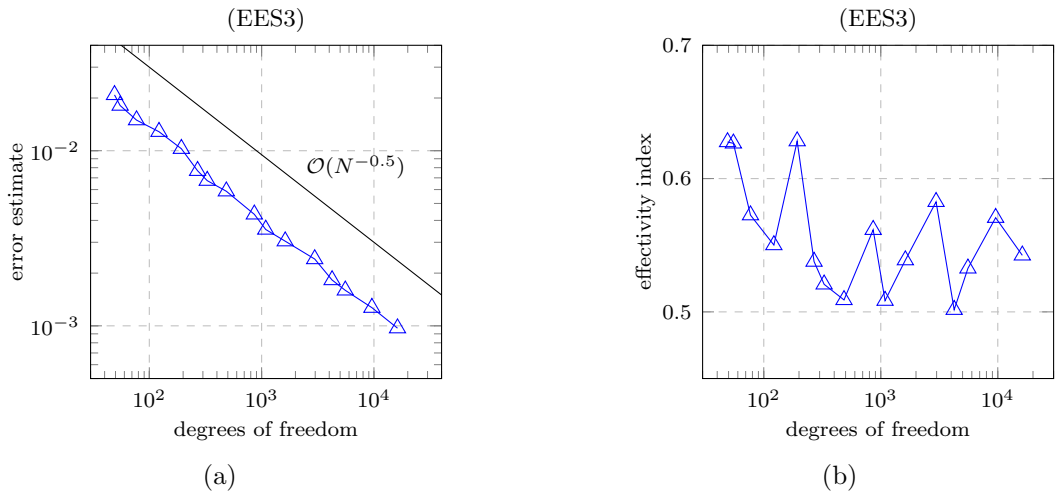

\rbl{In our second experiment,}
we ran the adaptive algorithm driven by two-level error estimates (i.e., using (EES3)) and
employed the {edge-based} D\"orfler marking~\eqref{marking:doerfler} with $\theta = 0.5$
(we use the same tolerance and the same coarse grid as before).
In this experiment, for the error estimate at each iteration of the adaptive algorithm,
we calculated the effectivity index (i.e., the ratio between the error estimate and 
\rbl{a surrogate approximation of the true error,
computed by running the adaptive algorithm with $P_2$ approximation
with a tighter tolerance of  {\tt tol = 2e-5}}).
In Figure~\ref{fig:ex1:exp2}(a), we plot the energy error estimates at each iteration.
Comparing this plot with the one in Figure~\ref{fig:ex1:exp1:converge}(c),
we can see improvements in terms of the monotonicity of the error decay and in terms of the number of iterations.
The \rbl{computed} effectivity indices for each iteration of the algorithm are plotted in Figure~\ref{fig:ex1:exp2}(b).

\rbl{In our third experiment, we repeated the previous adaptive procedure  with two smaller stopping
 tolerances. The overall computational times\footnote{\rbl{All timings
reported in this paper were recorded using
MATLAB \abrev{9.2 (R2017a)}} on a laptop with Intel Core i7 2.9GHz CPU and 16GB of RAM.}
 together with  the contributions for each \rbl{of the components in~\eqref{eq:modules}}  are 
reported in Table~\ref{tbl:ex1:exp2:times}}.

\smallskip

\begin{table}[t!]
\setlength\tabcolsep{16pt} 
\begin{center} 
\footnotesize{
\renewcommand{\arraystretch}{1.15}
\begin{tabular}{l   l  l   l  l    l  l   }
\toprule
\multicolumn{7}{c}{Adaptive FEM with \emph{edge-based} D\"orfler marking ($\theta=0.5$)}\\
\midrule
&\multicolumn{2}{l}{\tt tol=1e-3}	
&\multicolumn{2}{l}{\tt tol=5e-4}
&\multicolumn{2}{l}{\tt tol=1e-4}\\[-1pt]
\midrule 
$L$					&\multicolumn{2}{l}{$16$}
					&\multicolumn{2}{l}{$19$}	
					&\multicolumn{2}{l}{$27$}\\[-3pt]
$\eta_L$			&\multicolumn{2}{l}{$9.661$\text{e-}$04$}	
					&\multicolumn{2}{l}{$4.923$\text{e-}$04$}	
					&\multicolumn{2}{l}{$8.168$\text{e-}$05$}\\[-3pt]
$\abrev{n_L}$				&\multicolumn{2}{l}{16,261}	
					&\multicolumn{2}{l}{63,864}	
					&\multicolumn{2}{l}{2,181,895}\\[-1pt]					
\cmidrule(lr){1-7}
$t$ ({\sf SOLVE})	&\multicolumn{2}{l}{$0.167$}
					&\multicolumn{2}{l}{$0.811$}
					&\multicolumn{2}{l}{$49.567$}\\[-3pt]
$t$ ({\sf ESTIMATE})&\multicolumn{2}{l}{$0.324$}
					&\multicolumn{2}{l}{$2.013$}
					&\multicolumn{2}{l}{$107.107$}\\[-3pt]
$t$ ({\sf MARK})	&\multicolumn{2}{l}{$0.004$}
					&\multicolumn{2}{l}{$0.016$}
					&\multicolumn{2}{l}{$0.656$}\\[-3pt]
$t$ ({\sf REFINE})	&\multicolumn{2}{l}{$0.101$}
					&\multicolumn{2}{l}{$0.366$}
					&\multicolumn{2}{l}{$15.394$}\\[-1pt]
\cmidrule(lr){1-7}
$t$ (overall)		&\multicolumn{2}{l}{$3.090$}
					&\multicolumn{2}{l}{$7.937$}
					&\multicolumn{2}{l}{$439.862$}\\
\bottomrule 
\end{tabular}
\caption{
Example 1:
the outputs of running the adaptive algorithm employing the error estimation strategy~(EES3)
for three different tolerances.
Here, $L$, $\eta_L$, and \abrev{$n_L$} are as in Table~\ref{tbl:ex1:exp1:data}.
All times $t$ are in seconds and the timings for individual modules are \rbl{recorded at the final adaptive step}.
}
\label{tbl:ex1:exp2:times}
} 
\end{center}                                                                   
\end{table}

{\bf Example 2.} \
This \rbl{experiment} addresses the question posed by Nick Trefethen and Abi Gopal to the
 readers of the NA Digest in November 2018
(see~\cite{Trefethen_NA-Digest, GopalTrefethen19}).
The community was challenged 
to compute (to a high accuracy) the point-value close to singularity for a harmonic function in the L-shaped domain.
More precisely, let $D = (-1,1)^2 \setminus (-1,0]^2$ and consider the following  problem:
\be \label{eq:ex2:problem}
\begin{aligned}
           -\rbl{\nabla^2}  u(\bx) & = 0, \qquad && \bx = (x_1,x_2) \in D,\\[3pt]
           u(\bx) & = (1-x_1)^2, \qquad && \bx \in \partial D.
\end{aligned} 
\ee
The goal is to compute $u(0.01,0.01)$ to at least 8-digit accuracy
\rev{(the exact value is $1.02679192610\ldots$; see~\cite{GopalTrefethen19})}.

\begin{figure}[b!] 
\centering
\subfloat[][]
{
\includegraphics[scale = 0.95]{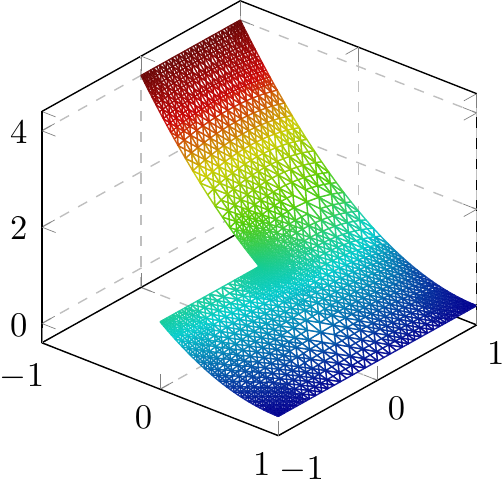}
}
\hfill   
\subfloat[][] 
{
\begin{tikzpicture}
\pgfplotstableread{data/exp2/eightdigit-conv.dat}{\first}
\begin{loglogaxis}
[
width = 6.0cm, height=5.5cm,
title={}, 
title style={font={\fontsize{9pt}{12pt}\selectfont},yshift=-1.0ex},
xlabel={degrees of freedom}, 
xlabel style={font={\fontsize{8pt}{12pt}\selectfont}},
ylabel={error estimate},
ylabel style={font={\fontsize{8pt}{12pt}\selectfont},yshift=-1.0ex},
ymajorgrids=true, xmajorgrids=true, grid style=dashed, 
xmin = (1)*10^(2), 	xmax = (6)*10^(5),
ymin = (2)*10^(-5),	ymax = (1.5)*10^(-1),
font={\fontsize{8pt}{12pt}\selectfont}, 
legend style={legend pos=south west, legend cell align=left, fill=none, draw=none, font={\fontsize{9pt}{12pt}\selectfont}}
]
\addplot[blue,mark=triangle,mark size=3.5pt]	table[x=dofs, y=error]{\first};
\addplot[black,solid,domain=10^(1):10^(7)]    { (65.5)*x^(-1) };
\node at (axis cs:1e4,6e-3)	[anchor=south west]{$\mathcal{O}(N^{-1})$};
\end{loglogaxis}
\end{tikzpicture}
} 
\hfill 
\subfloat[][]
{
\includegraphics[scale = 1.25]{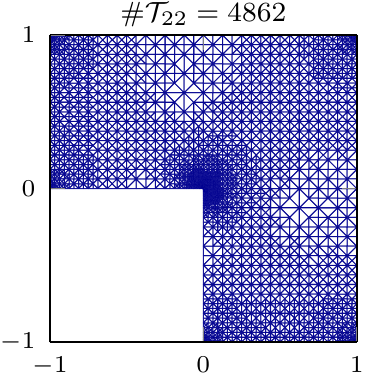}
}
\caption{
Example 2: 
(a) Galerkin solution to problem \eqref{eq:ex2:problem};
(b) error estimates at each iteration of the adaptive algorithm with $P_2$ approximations 
and {\tt tol = 4e-5};
(c) adaptively refined mesh generated by the algorithm for \rbl{an intermediate tolerance of~{\tt 1e-3}}.
}
\label{fig:ex2}
\end{figure}
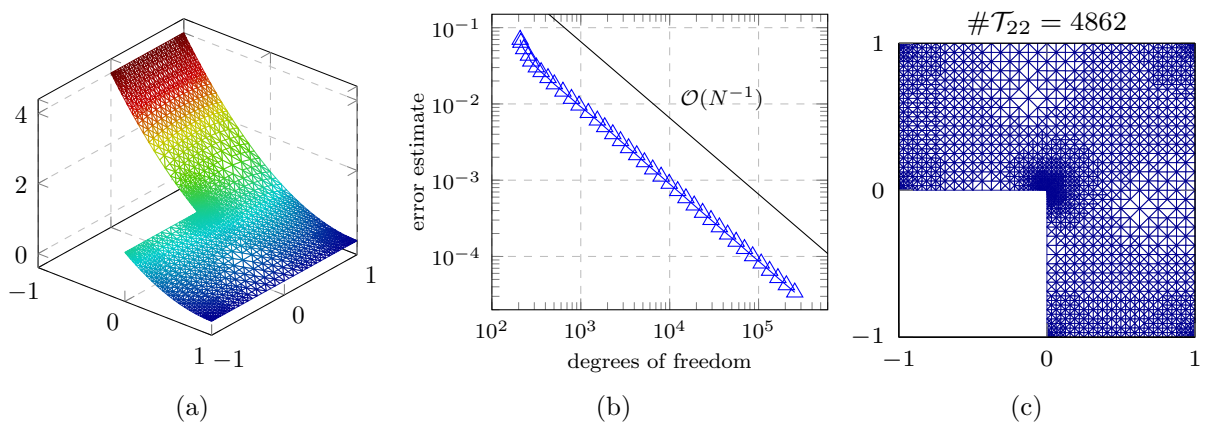

In this example, we set the stopping tolerance   {\tt tol = 4e-5} and
\rbl{ran} the adaptive FEM algorithm \rbl{with $P_2$ approximation \abrev{together with element-based 
D\"orfler marking with}} $\theta = 0.5$.
The \rbl{prescribed tolerance was satisfied} after 38 iterations (final number of d.o.f. was 253,231, run time  59.6 sec),
giving the value $u(0.01,0.01) \approx 1.02679192311$, \rbl{which} is accurate to 9 digits.
Figure~\ref{fig:ex2} depicts the finite element solution to problem~\eqref{eq:ex2:problem}
and shows the convergence plot for the estimated energy errors (\rbl{together}  with the optimal rate) 
as well as the mesh refinement pattern, plotted here for an intermediate tolerance.
%

\section{Goal-oriented adaptivity} \label{sec:goafem}
\setcounter{equation}{0}

The error estimation strategies described in the previous section provide a mechanism
for controlling approximation errors in the (global) energy norm.
However, in many practical applications, simulations often target a specific quantity of interest
\rbl{(typically, a local feature of the solution) called the {\it goal functional}. 
In such cases, the energy norm of the error is likely to be of limited interest.
The implementation of  {\it goal-oriented} error estimation and adaptivity in T-IFISS is the focus of this section.
We also discuss a representative case \abrev{study} to show the efficiency of the adopted approach.} 

\subsection{\rev{Goal-oriented error estimation in the abstract setting}} \label{subsec:goafem:estimate}

We start by describing a general idea of the goal-oriented error estimation strategy implemented in T-IFISS.
Let $V$ be a Hilbert space and denote by $V'$ its dual space.
Let $B: V \times V \to \RR$ be a continuous, elliptic, and symmetric bilinear form 
with the associated energy norm $\bnorm{\cdot}$, i.e., $\bnorm{v}^2 := B(v,v)$ for all $v \in V$.
Given two continuous linear functionals $F, G \in V'$,
our aim is to approximate $G(u)$, where $u \in V$ is the unique solution
of the \emph{primal problem}:
\begin{equation} \label{eq:primal:problem}
      B(u,v) = F(v)
\quad \text{for all $v \in V.$}
\end{equation}
To this end, the standard approach (see, e.g., \cite{MR1352472, prudhommeOden1999, beckerRannacher2001, gilesSuli2002})
considers $z \in V$ as a unique solution to the \emph{dual problem}:
\begin{equation} \label{eq:dual:problem}
      B(v,z) = G(v) \quad \text{for all $v \in V$.}
\end{equation}
Let $V_h$ be a finite dimensional subspace of $V$. Let 
$u_h \in V_h$ (resp., $z_h \in V_h$) be a unique Galerkin approximation
of the solution to the primal (resp., dual) problem, i.e., 
\[
     B(u_h,v_h) = F(v_h) \quad
     (\text{resp., } B(v_h,z_h) = G(v_h)) \qquad \text{for all $v_h \in V_h.$}
\]
Then, it follows that
\begin{equation} \label{eq:errorInequality}
   \lvert G(u)-G(u_h) \rvert
   = \lvert B(u-u_h,z) \rvert
   = \lvert B(u-u_h,z-z_h)\rvert
       \leq \bnorm{u-u_h}\,\bnorm{z-z_h},
\end{equation}
where the second equality holds due to Galerkin orthogonality. 

Assume that $\mu = \mu(u_h)$ and $\zeta = \zeta(z_h)$ are reliable estimates
for the energy errors $\bnorm{u-u_h}$ and $\bnorm{z-z_h}$, respectively, i.e.,  
\begin{equation*} 
      \bnorm{u-u_h} \lesssim \mu
      \quad\text{and}\quad
      \bnorm{z-z_h} \lesssim \zeta
\end{equation*}
(here, $a \lesssim b$ means the existence of a generic positive constant $C$ such that $a \le C b$).
Hence, inequality \eqref{eq:errorInequality} implies that the product $\mu \, \zeta$ is a reliable error
estimate for the approximation error in the goal functional:
\begin{equation} \label{goal:error:reliability}
      \lvert G(u)-G(u_h) \rvert  
      \lesssim \mu \, \zeta.
\end{equation}

\subsection{\rev{Marking in goal-oriented adaptivity}} \label{subsec:goafem:mark}

Having computed two Galerkin solutions $u_h$, $z_h$, the corresponding energy error estimates
$\mu(u_h)$, $\zeta(z_h)$, and the reliable estimate $\mu(u_h) \, \zeta(z_h)$ 
of the error in the goal functional (see~\eqref{goal:error:reliability}),
a goal-oriented adaptive FEM (GOAFEM) algorithm proceeds by executing the 
$\mathsf{MARK}$ and $\mathsf{REFINE}$ modules
of the standard adaptive loop~\eqref{eq:modules}.

While the module $\mathsf{REFINE}$ simply performs local mesh refinement 
\rbl{as explained in the previous section,}
the marking procedure in the GOAFEM algorithm requires \rbl{special care}.
\rbl{Specifically,} since the error in the goal functional is controlled by the product of the energy error estimates
for two Galerkin approximations (the primal and dual ones),
the edge-marking for bisection (or, the element-marking for refinement)
must take into account the local error indicators associated \abrev{with} {\it both} approximations.
\rbl{Thus,}  given the set of local error indicators associated 
with the primal (resp., dual) Galerkin solution $u_h$ (resp., $z_h$),
let $\mathcal{M}^u$ (resp., $\mathcal{M}^z$) denote the set of element edges that would be marked for bisection
in order to enhance this Galerkin solution
(to that end, one can use, e.g., the D\"orfler marking strategy, see~\eqref{marking:doerfler}).
There exist several strategies \rbl{for combining} the two 
sets $\mathcal{M}^u$ and $\mathcal{M}^z$ into a single marking set $\mathcal{M}$
that is used for mesh refinement in the goal-oriented adaptive algorithm.
Four \rbl{different strategies} are implemented in T-IFISS:

(GO--MARK1)
following~\cite{hp16},
the marking set $\mathcal{M}$ is simply the union of $\mathcal{M}^u$ and $\mathcal{M}^z$;

(GO--MARK2) 
following~\cite{ms09},
the marking set $\mathcal{M}$ is defined as the set of minimal cardinality between $\mathcal{M}^u$ and $\mathcal{M}^z$;

(GO--MARK3)
following~\cite{bet11}, the set $\mathcal{M}$ is obtained by performing D\"orfler marking on the set 
of combined error indicators $\beta(E)$ associated with edges of the triangulation,~where
\[
   \beta(E) := \big(\mu^2_E \, \zeta^2 + \zeta^2_E \, \mu^2\big)^{1/2}
\]
and $\mu_E$ (resp., $\zeta_E$) is the local contribution to $\mu$ (resp., $\zeta$)
associated with the edge $E$;

(GO--MARK4)
this marking strategy is a modification of (GO--MARK2);
following~\cite{fpv16}, we compare the cardinality of $\M^u$ and that of $\M^z$ to define
\begin{align}
        \M_\star := \M^u \quad \hbox{and} \quad \M^\star := \M^z \qquad &
        \hbox{if $\#\M^u \le \#\M^z$},
        \nonumber
		\\ 
        \M_\star := \M^z \quad \hbox{and} \quad \M^\star := \M^u \qquad &
        \hbox{otherwise};
        \nonumber
\end{align}
the marking set $\M$ 
is then defined as the union of $\M_\star$ and 
those $\#\M_\star$ edges of $\M^\star$ that have the largest error indicators.

\rbl{Comparing these four strategies,}
it is proved in~\cite[Theorem~13]{fpv16} that the GOAFEM algorithm employing marking 
strategies (GO--MARK2)--(GO--MARK4)
generates approximations that converge with {\it optimal} algebraic rates,
whereas only {\it suboptimal} convergence rates have been proved for marking strategy (GO--MARK1);
cf. \cite[Remark~4]{fpv16} and \cite[Section~4]{hp16}.
\rbl{The numerical results  in~\cite{fpv16} suggest that (GO--MARK4) is 
more effective than the original strategy (GO--MARK2) in terms of the overall computational~cost.
Our own experience is that  (GO--MARK4) is a competitive strategy
in every example that has been tested. Consequently, we have made it the default option within the code.
}

\subsection{\rev{Numerical case study}} \label{subsec:goafem:numerics}

\rev{In order to demonstrate the effectiveness of the goal-oriented adaptive strategy
described in the previous subsection, let us consider the following test example.
}


\smallskip

\begin{figure}[t!] 
\centering
\subfloat[][]
{
\includegraphics[scale=0.84]{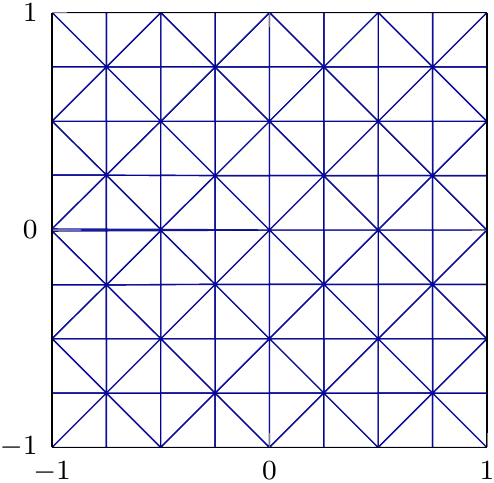}} 
\hspace{1.2cm}
\subfloat[][] 
{
\includegraphics[scale=0.85]{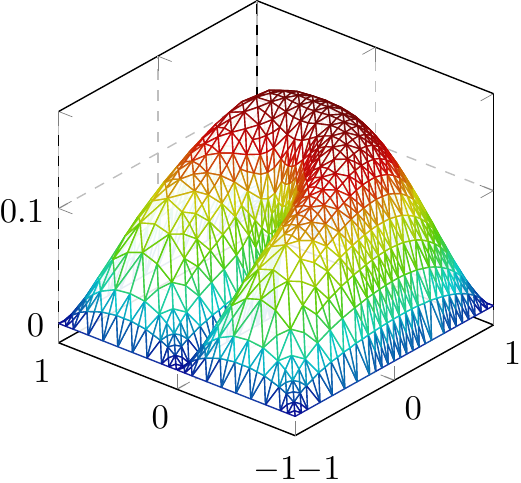}} 
\caption{
Example 3: (a) initial coarse mesh 
in the GOAFEM algorithm; (b) the primal Galerkin~solution.
}
\label{fig:ex4:primal:sol}
\end{figure}

{\bf Example~3.}\  \rbl{Let us} consider the model problem given by~\eqref{eq:ex1:problem} with $A = \text{Id}$
on \rbl{a slit domain
$D_\delta = (-1,1)^2 \setminus \overline{T}_\delta$,
where $T_\delta = \text{conv}\{(0,0), (-1,\delta), (-1,-\delta)\}$ with $\delta = 0.005$.}
%
It is known that solution $u$ to the (primal) problem in this example exhibits a singularity
induced by the slit in the domain (see Figure~\ref{fig:ex4:primal:sol}(b)).
Our aim, however, is to demonstrate the capability of the software to approximate the value of $u$
at some fixed point $\bx_0 \in D$ away from the slit
(in the experiments below, we set $\bx_0 = (0.4,-0.5)$).
In order to define the corresponding bounded goal functional $G$, it is common to
fix a sufficiently small $r >0$ and first introduce the \emph{mollifier} $g_0$ 
as follows (cf.~\cite{prudhommeOden1999}):
\begin{equation} \label{mollifier}
   g_0(\bx) = g_0(\bx; \bx_0, r) := 
   \begin{cases}
      C \exp\left(- \frac{r^2}{r^2 - |\bx - \bx_0|^2} \right)	& \text{if $|\bx - \bx_0| < r$}, \\
      0 								& \text{otherwise.}
   \end{cases}
\end{equation}
Here, 
the constant $C$ is chosen such that
$
   \int_D g_0(\bx)\, \dx = 1
$
(for sufficiently small $r$ such that $\supp(g_0(\bx;\bx_0,r)) \subset D$, one has
$C \approx 2.1436\,r^{-2}$; see, e.g., \cite{prudhommeOden1999}).
%
Then, the functional $G$ in~\eqref{eq:dual:problem} reads as
\begin{equation*}
  G(v) = \int_D g_0(\bx) v(\bx) \, \dx \quad \text{for all } v \in H^1_0(D). 
\end{equation*}
Note that if $u(\bx)$ is continuous in a neighborhood of $\bx_0$, then
$G(u)$ converges to the point value $u(\bx_0)$ as $r$ tends to zero.

\begin{figure}[b!]
\centering
\footnotesize
\includegraphics[scale=0.84]{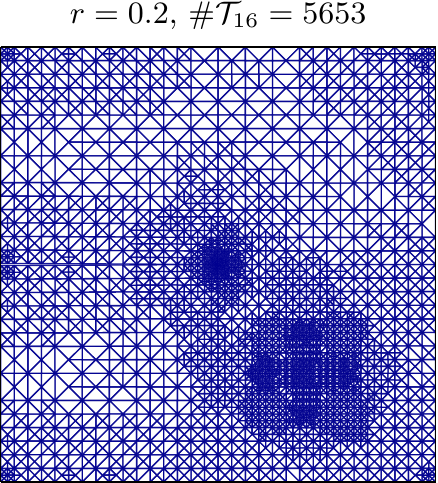} 				\hspace{1.15cm}
\includegraphics[scale=0.84]{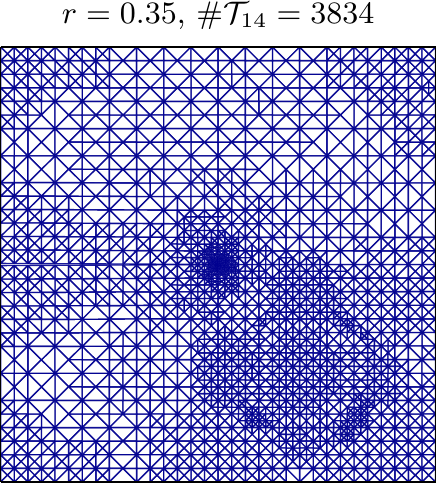} 				\hspace{1.15cm}
\includegraphics[scale=0.84]{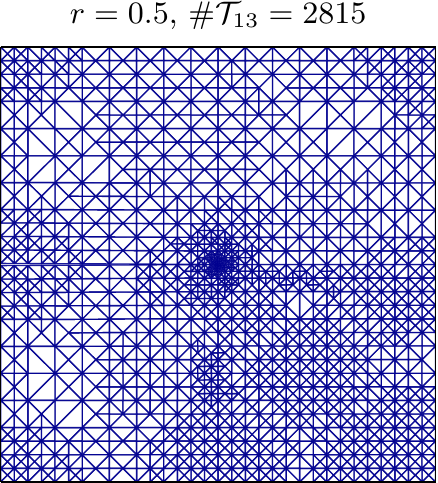}   				\\
\vspace{6pt}
\includegraphics[scale=0.8]{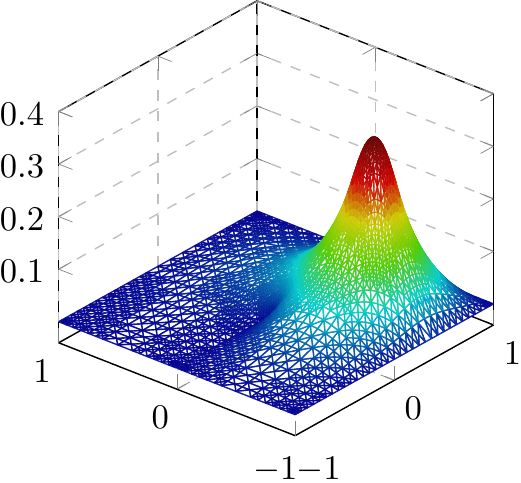} 	\hspace{0.75cm}
\includegraphics[scale=0.8]{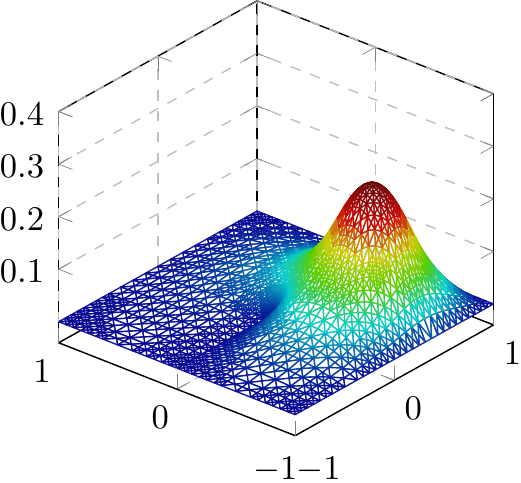} 	\hspace{0.7cm}
\includegraphics[scale=0.8]{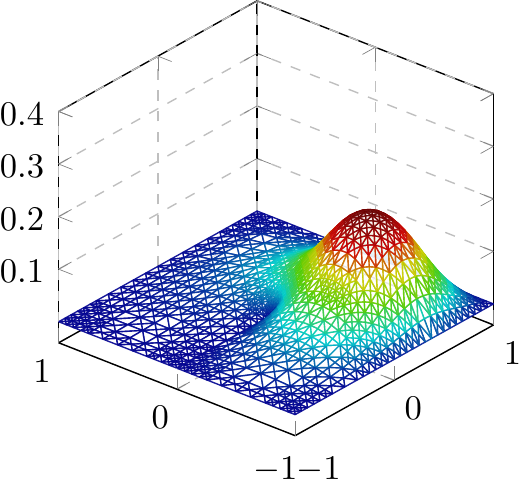}
\caption{
\rbl{Example~3}: adaptively refined triangulations (top row) and the dual Galerkin solutions (bottom row)
computed using the mollifier $g_0$ in~\eqref{mollifier} with $r = 0.2,\, 0.35,\, 0.5$.
}
\label{fig:ex4:sequence:meshes}
\end{figure}

\rbl{We started with the coarse triangulation 
depicted in Figure~\ref{fig:ex4:primal:sol}(a) in all the experiments.}
\rbl{In the first experiment, we fixed the stopping tolerance  to be  {\tt tol = 3e-4}  and ran} the GOAFEM algorithm to 
compute dual Galerkin solutions for different values of the radius $r$ in~\eqref{mollifier}.
\rbl{For the $\mathsf{SOLVE}$ step, we used $P_1$ approximations for both primal and dual solutions.
Within the $\mathsf{ESTIMATE}$ module, the energy errors in both solutions were 
estimated using the two-level error estimation strategy (EES3) described earlier.
Given the error indicators for primal and dual solutions, the algorithm employed the edge-based
D\"orfler marking~\eqref{marking:doerfler} with $\theta = 0.3$ in combination with 
the strategy (GO--MARK4)  above.}
Figure~\ref{fig:ex4:sequence:meshes} shows the refined triangulations (top row) and the corresponding
dual Galerkin solutions (bottom row) for $r = 0.2,\, 0.35,\, 0.5$.
We note that the triangulations generated by the algorithm adapt to the features of both primal and dual solutions:
the triangulations are refined in the vicinity of each corner
(with particularly strong refinement near the tip of the slit) and
in a neighborhood of~$\bx_0$ (with stronger refinement for smaller values of $r$).

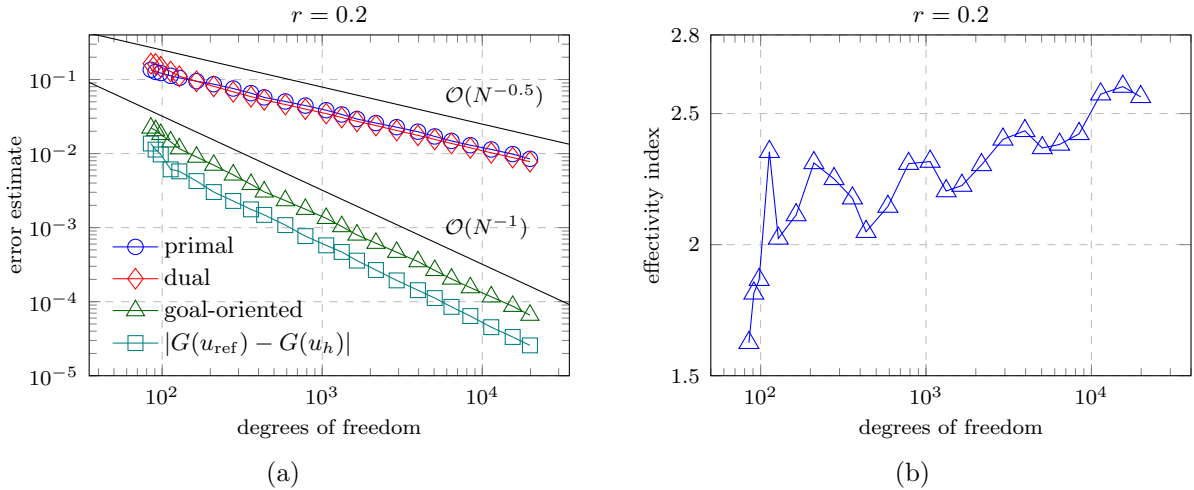
\begin{figure}[t!] 
\centering
\subfloat[][] 
{
\begin{tikzpicture}
\pgfplotstableread{data/exp4/r=02.dat}{\first}
\begin{loglogaxis}
[
width = 7.9cm, height=6.1cm,
title={$r=0.2$}, 
title style={font={\fontsize{9pt}{12pt}\selectfont},yshift=-1.0ex},
xlabel={degrees of freedom}, 
xlabel style={font={\fontsize{8pt}{12pt}\selectfont}},
ylabel={error estimate},
ylabel style={font={\fontsize{8pt}{12pt}\selectfont},yshift=-1.0ex},
ymajorgrids=true, xmajorgrids=true, grid style=dashed, 
xmin = (3.5)*10^(1), 	xmax = (3.5)*10^(4),
ymin = (1)*10^(-5),		ymax = (4)*10^(-1),
font={\fontsize{8pt}{12pt}\selectfont}, 
legend style={legend pos=south west, legend cell align=left, fill=none, draw=none, font={\fontsize{9pt}{12pt}\selectfont}}
]
\addplot[blue,mark=o,mark size=2.9pt]				table[x=dofs, y=errprimal]{\first};
\addplot[Red,mark=diamond,mark size=3.9pt]			table[x=dofs, y=errdual]{\first};
\addplot[darkGreen,mark=triangle,mark size=3.9pt]	table[x=dofs, y=errprod]{\first};
\addplot[teal,mark=square,mark size=2.8pt]	table[x=dofs, y=truegerr]{\first};
\addplot[black,solid,domain=10^(1):10^(7)]    { (2.5)*x^(-0.5) };
\node at (axis cs:5e3,3e-2)	[anchor=south west]{$\mathcal{O}(N^{-0.5})$};
\addplot[black,solid,domain=10^(1):10^(7)]    { (3.2)*x^(-1) };
\node at (axis cs:5e3,5e-4)	[anchor=south west]{$\mathcal{O}(N^{-1})$};
\legend{
{primal},
{dual},
{goal-oriented},
{$|G(u_{\rm ref}) - G(u_h)|$}
}
\end{loglogaxis}
\end{tikzpicture}
} 
\hfill 
\subfloat[][]
{
\begin{tikzpicture}
\pgfplotstableread{data/exp4/r=02.dat}{\first}
\begin{semilogxaxis}
[
width = 7.9cm, height=6.1cm,
title={$r=0.2$}, 
title style={font={\fontsize{9pt}{12pt}\selectfont},yshift=-1.0ex},
xlabel={degrees of freedom}, 
xlabel style={font={\fontsize{8pt}{12pt}\selectfont}},
ylabel={effectivity index},
ylabel style={font={\fontsize{8pt}{12pt}\selectfont},yshift=-1.0ex},
ymajorgrids=true, xmajorgrids=true, grid style=dashed, 
xmin = (5)*10^(1), 	xmax = (4)*10^4,
ymin = 1.5,			ymax = 2.8,
ytick={1.5, 2, 2.5, 2.8},
font={\fontsize{8pt}{12pt}\selectfont}, 
legend style={legend pos=south west, legend cell align=left, fill=none, draw=none, font={\fontsize{9pt}{12pt}\selectfont}}
]
\addplot[blue,mark=triangle,mark size=4.3pt]	table[x=dofs, y=effindices]{\first};
\end{semilogxaxis}
\end{tikzpicture}
}
\caption{
\rbl{Example 3}:
(a) error estimates at each iteration of the GOAFEM algorithm employing the marking strategy (GO--MARK4);
(b) the associated effectivity indices. 
}
\label{fig:ex4:converge:effectivity}
\end{figure}

Focusing now on the case $r=0.2$, 
\rbl{we set  {\tt tol = 8e-5} in the second experiment and ran the GOAFEM  algorithm without changing the 
settings.} 
\rbl{The results we obtained are shown  in Figure~\ref{fig:ex4:converge:effectivity}.}
In Figure~\ref{fig:ex4:converge:effectivity}(a), we plot the energy error 
estimates for primal and dual Galerkin approximations, 
the estimates of the error in the goal functional, as well as the reference errors in the goal functional
(i.e., $|G(u_{\rm ref}) - G(u_h)|$) at each iteration of the GOAFEM algorithm.
Here, the reference Galerkin solution $u_{\rm ref}$ is computed using the triangulation 
obtained by two uniform refinements of the final triangulation generated by the GOAFEM algorithm.
We observe that all error estimates as well as the reference error in the goal functional converge with optimal rates.
The effectivity indices for the goal-oriented error estimation at each iteration of the algorithm are plotted 
in Figure~\ref{fig:ex4:converge:effectivity}(b).
This plot shows that the product of energy error estimates for the primal and dual Galerkin solutions 
provides a \rbl{reasonably} accurate
estimate for the error in approximating the goal functional~$G(u)$.

\input{section3.tex}


\section{Extensions and future developments} \label{sec:extensions}

\rbl{
The T-IFISS software framework provides many opportunities for experimentation and exploration.
It is also an invaluable teaching tool for numerical analysis and 
computational engineering courses with an emphasis on contemporary finite element analysis.}
Stochastic T-IFISS has been recently extended to incorporate
the goal-oriented error estimation and adaptivity 
in the context of stochastic Galerkin approximations for parametric elliptic PDEs
(see~\cite{bprr19} for details of the algorithm and numerical results).
The toolbox has also been used for numerical testing of the adaptive algorithm proposed
for parameter-dependent linear elasticity problems; 
see~\cite{KhanPS_19_RPS, KhanBespalovPowellSilvester}.
Future developments would include the extension to problems with non-affine parametric 
representations of inputs (see~\cite{BespalovX_AEE})
and the implementation of the multilevel adaptive SGFEM algorithm
 in the spirit of~\cite{egsz14, CrowderPB_19_EAM}.

%

\bigskip

\noindent{\bf Acknowledgement.}
The authors are grateful to Qifeng Liao (ShanghaiTech University) for his inputs to T-IFISS project at its early stages.
The authors would also like to thank Dirk Praetorius and Michele Ruggeri (both at Technical University of Vienna)
for their contributions to the development of the toolbox components for goal-oriented error estimation and adaptivity.

\newpage

\bibliographystyle{siam}
\bibliography{references}

\end{document}

%% file: section3.tex

\section{\rbl{Adaptive FEM for parametric PDEs}} \label{sec:asgfem}
\setcounter{equation}{0}

\newcommand{\dpi}{\mathrm{d}\pi}
\newcommand{\XX}{\mathbb{X}}
\newcommand{\YY}{\mathbb{Y}}


\rev{In this section, we address the design of components of the adaptive algorithm when working in a parametric PDE setting.}
Before discussing the details of \rbl{our implementation, 
we formulate} the model problem with parametric input data
and briefly describe the idea of stochastic Galerkin FEM (SGFEM).
\rbl{Readers interested in theoretical aspects of SGFEM are referred  to~\cite{gs91}, \cite{dbo01} and~\cite{ btz04}.}

\subsection{\rbl{Stochastic Galerkin FEM}} \label{sec:asgfem:model}

Let $D \subset \RR^2$ be a Lipschitz domain (called the physical domain) with polygonal boundary $\partial D$
and let $\G := \prod_{m=1}^\infty [-1,1]$ denote the infinitely-dimensional hypercube (called the parameter domain).
We consider the elliptic boundary value problem
\begin{equation} \label{eq:strongform}
\begin{aligned}
 -\nabla \cdot (a \nabla u) &= f \quad &&\text{in } D \times \G,\\
 u & = 0 \quad &&\text{on } \partial D \times \G,
\end{aligned}
\end{equation}
where the scalar coefficient $a$ 
(and, hence, the solution~$u$)
depends on a countably infinite number of scalar parameters, i.e.,
$a = a(\bx, \boldy)$ 
and $u = u(\bx, \boldy)$ with $\bx \in D$, $\boldy \in \G$,
and the differentiation in $\nabla$ is with respect to $\bx = (x_1,x_2)$.
We assume that $f = f(\bx) \in H^{-1}(D)$
and that the parametric coefficient $a$ has affine dependence on the parameters, i.e.,
\begin{equation} \label{eq1:a}
   a(\bx,\boldy) = a_0(\bx) + \sum_{m=1}^\infty y_m a_m(\bx)
   \quad \text{for } \bx \in D \text{ and } \boldy = (y_m)_{m\in\NN} \in \G,
\end{equation}
%
where the scalar functions $a_m \in W^{1,\infty}(D)$ ($m \in \NN_0$) 
satisfy the following~inequalities:
\begin{equation}\label{eq2:a}
 0 < a_0^{\min} \le a_0(\bx) \le a_0^{\max} < \infty
 \quad \text{for almost all } \bx \in D
\end{equation}
and
\begin{equation}\label{eq3:a}
 \tau := \frac{1}{a_0^{\rm min}} \, \sum_{m=1}^\infty \|a_m\|_{L^\infty(D)} < 1.
\end{equation}
The weak formulation of~\eqref{eq:strongform} is posed in the framework of the Bochner space $V := L^2_\pi(\G;H^1_0(D))$.
Here, $\pi \,{=}\, \pi(\boldy)$ is a probability measure on $(\G,\mathcal{B}(\G))$ with
$\mathcal{B}(\G)$ being the Borel $\sigma$-algebra on $\G$, and
we assume that $\pi(\boldy)$ is the product of symmetric Borel probability measures $\pi_m$ on~$[-1,1]$,
i.e., $\pi(\boldy) = \prod_{m=1}^\infty \pi_m(y_m)$.
For a given $f \in H^{-1}(D)$, the weak solution $u \in V$ satisfies
\begin{equation} \label{eq:weakform}
B(u,v) = F(v) := \int_\G \int_D f(\bx) v(\bx,\boldy) \, \dx \, \dpi(\boldy)
\quad \text{for all } v \in V.
\end{equation}
Here,
\begin{align}
   B(u,v) &:= B_0(u,v) + \sum_{m=1}^\infty \int_\G \int_D y_m a_m(\bx) \nabla u(\bx,\boldy)\cdot\nabla v(\bx,\boldy) \, \dx \, \dpi(\boldy),
   \nonumber
\\
   B_0(u,v) &:= \int_\G \int_D a_0(\bx) \nabla u(\bx,\boldy)\cdot\nabla v(\bx,\boldy) \, \dx \, \dpi(\boldy).
   \label{eq:B0}
\end{align}

The following three coefficient expansions (CEs) of the type~\eqref{eq1:a} are implemented in Stochastic T-IFISS.

(CE1)
Let $D=(-1,1)^2$ and
suppose that the coefficient $a = a(\bx,\boldy)$ in~\eqref{eq:strongform} is a parametric representation
of a second-order random field with prescribed mean $\EE[a]$ and covariance function~$\Cov[a]$.
Assume that $\Cov[a]$ is the separable exponential covariance function given by 
\begin{equation*} 
       \Cov[a](\bx,\bx') = \sigma^2 \exp\left(-\frac{|x_1 - x_1'|}{l_1} - \frac{|x_2 - x_2'|}{l_2}\right), 
\end{equation*}
where $\bx,\, \bx' \in D$, $\sigma$ denotes the standard deviation, 
and $l_1, l_2$ are correlation~lengths.
In this case, $a(\bx,\boldy)$ can be written in the form~\eqref{eq1:a}
using the Karhunen--\rev{Lo{\`e}ve}-type expansion~\cite{lps2014}
such that
\begin{equation*} 
   a_0(\bx) := \EE[a](\bx),\quad
   a_m(\bx) := c \, \sqrt{\lambda_m}\,\varphi_m(\bx),\quad m\in \NN,
\end{equation*}
where $\{(\lambda_m, \varphi_m)\}_{m=1}^\infty$ are the eigenpairs of the 
operator $\int_D \Cov[a](\bx,\bx') \varphi(\bx') \dx'$, 
and the constant $c > 0$ 
is chosen such that ${\rm Var}(c \, y_m)=1$ for all $m \in \NN$.
Note that analytical expressions for $\lambda_m$ and $\varphi_m$ for the square domain $D$
follow from the corresponding formulas derived in~\cite[pp.~28--29]{gs91} for the one-dimensional case.

(CE2)
Following~\cite[Section~11.1]{egsz14}, we set $a_0(\bx) := 1$, $\bx \in D$, and choose the coefficients $a_m(\bx)$ in~\eqref{eq1:a} 
to represent planar Fourier modes of increasing total order:
\begin{equation} \label{Eigel:coeff}
   a_m(\bx) := \alpha_m \cos(2\pi\beta_1(m) \, x_1) \cos(2\pi\beta_2(m)\, x_2)
\quad \text{for all $m \in \NN$},
\end{equation}
where $\alpha_m := A m^{-\tilde\sigma}$ are the amplitudes of the coefficients, $\tilde\sigma > 1$,
the constant $A$ is chosen such that $\tau = A \zeta(\tilde\sigma) = 0.9$ (here, $\zeta$ denotes the Riemann zeta function),
and $\beta_1$, $\beta_2$ are defined as
\[
\beta_1(m) := m - k(m)(k(m) + 1)/2\ \ \hbox{and}\ \
\beta_2(m) := k(m) - \beta_1(m),
\]
with $k(m) := \lfloor -1/2  + \sqrt{1/4+2m}\rfloor$ for all $m \in \NN$.
%
The following two values of the decay parameter $\tilde\sigma$ are used in test problems:
$\tilde\sigma = 2$ (slow decay) and $\tilde\sigma = 4$ (fast decay).

(CE3)
Finally, we consider the following parametric coefficient (cf.~\cite[Example 9.37]{lps2014}):
\be \label{PowellRandomField_1}
	a(\bx, \tilde\boldy) = 1 + 
	c\, \sum_{i=0}^\infty \sum_{j=0}^{\infty}  \sqrt{\lambda_{ij} } \varphi_{ij}(\bx) y_{ij},\quad
	\tilde\boldy = (y_{ij})_{i,j \in \NN_0},
\ee
where
$\lambda_{ij} = \bar\lambda_i\, \bar\lambda_j$,
$\varphi_{ij}(\bx) = \bar\varphi_{i}(x_1) \, \bar\varphi_{j}(x_2)$,
$y_{ij} \in [-1,1]$ ($i,j \in \NN_0$) with
$\bar\lambda_0 := 1/2$, 
$\bar\varphi_{0}(t) := 1$,
$\bar\lambda_{k} := \frac{1}{2} \exp(-\pi k^2 \ell^2)$,
$\bar\varphi_{k}(t) := \sqrt{2} \cos(\pi k t)$ ($k \in \NN)$, and
the constant $c > 0$ 
is chosen such that ${\rm Var}(c \, y_{i,j})=1$ for all $i,j \in \NN_0$.

As shown in~\cite[Example 9.37]{lps2014},
the parametric representation~\eqref{PowellRandomField_1} stems
from the Karhunen-Lo\`eve expansion of a random field with the mean $\EE[a] = 1$ and
covariance function close to the isotropic covariance
$c(\xx) = (4\ell^2)^{-1}\,\exp(-\pi(x_1^2+x_2^2)/(4\ell^2))$,
where $\ell$ is the correlation length.
For our implementation, we set $\ell = 1$ and reorder the terms in the double sum in~\eqref{PowellRandomField_1}
to write the coefficient $a(\bx, \tilde\boldy)$ in the form~\eqref{eq1:a} with
\[ 
          a_0(\bx) := 1,\quad a_m(\bx) := c\, \sqrt{\lambda_m } \, \varphi_m(\bx),\quad y_m \in [-1,1],\quad m \in \NN,
\]
where $\lambda_1 \ge \lambda_2 \ge \ldots$.  
%
%

In each coefficient expansion (CE1)--(CE3), parameters $y_m$ ($m\in \NN$) are the images of
independent mean-zero random variables on $\Gamma_m = [-1,1]$.
The following two types of bounded random variables (RVs) are implemented in Stochastic T-IFISS.

(RV1) Uniformly distributed random variables.
In this case, $\dpi_m = \mathrm{d} y_m/2$ 
and the orthonormal polynomial basis in $L^2_{\pi_m}(\Gamma_m)$ is comprised of scaled Legendre~polynomials.

(RV2) Truncated Gaussian random variables.
In this case, $\dpi_m = p(y_m) \mathrm{d} y_m$ with
\begin{equation} 
   p(y_m) = \frac{\exp(-y_m^2 / (2\sigma_0^2))}{\sigma_0 \sqrt{2\pi} (2\Phi(1/\sigma_0) - 1)},\quad
   m \in \NN,
\end{equation}
where $\Phi(\cdot)$ is the Gaussian cumulative distribution function and
$\sigma_0$ is a parameter measuring the standard deviation.
The corresponding orthonormal polynomial basis in $L^2_{\pi_m}(\Gamma_m)$ is formed
by the so-called Rys polynomials; see, e.g., \cite[Example~1.11]{gautschi2004}.

A key observation that motivates the stochastic Galerkin FEM is that
the Bochner space $V = L^2_\pi(\G;H^1_0(D))$ is isometrically isomorphic to $H^1_0(D) \otimes L^2_\pi(\G)$.
Mimicking this tensor-product construction, the finite-dimensional subspace $V_{\XX\PP} \subset V$ is defined
as $V_{\XX\PP} := \XX \otimes \PP$;
here,
$\XX = \XX_h$ is a finite element space associated with a conforming triangulation $\mesh_h$ of $D$ and
$\PP = \PP_{\gotP}$ is a polynomial space on $\G$ associated with a finite index set $\gotP$.
Specifically,
\be \label{eq:fe:space}
   \XX = \XX_h :=  \Span\{\phi_i;\; i = 1,\ldots N_\XX\} \subset H^1_0(D),\quad
   N_\XX := \dim(\XX)
\ee
and
\[   
   \PP = \PP_{\gotP} :=
   \Span\Big\{ \sfP_\nu(\boldy) = \prod_{m\in\NN} \sfP_{\nu_m}^m(y_m);\;
   \nu \in \gotP
   \Big\} \subset L^2_\pi(\G)\ \ \text{with $\gotP \subset \gotI$},
\]   
where
$\{\sfP_n^m : n\in\NN_0 \}$ is an orthonormal basis of $L^2_{\pi_m}(-1,1)$ and
$\gotI$ denotes the countable set of finitely supported multi-indices, i.e.,
\[
   \gotI := \big\{\nu = (\nu_m)_{m\in\NN};\; \nu_m\in\NN_0 \text{ for all } m\in\NN,\ \#\supp(\nu) < \infty \big\}
\]
with $\supp(\nu):= \{m\in\NN;\; \nu_m\neq 0 \}$.

The Galerkin discretization of~\eqref{eq:weakform} reads as follows:
find $u_{\XX\PP} \in V_{\XX\PP}$ such that
\begin{equation}\label{eq:discrete_formulation}
   B(u_{\XX\PP}, v) = F(v) 
   \quad \text{for all } v \in V_{\XX\PP}.
\end{equation}
Note that $\dim(V_{\XX\PP}) = \dim(\XX_h) \rev{\times} \dim(\PP_\gotP)$.
Therefore, if a large number of random variables is used to represent the input data,
then computing high-fidelity stochastic Galerkin approximations
with standard polynomial subspaces on $\G$ (e.g., the spaces of tensor-product or complete polynomials)
becomes prohibitively expensive.
This motivates the development of adaptive SGFEM algorithms that incrementally refine
spatial ($\XX$-) and parametric ($\PP$-) components of Galerkin approximations
by iterating the standard adaptive loop~\eqref{eq:modules}.
%
\rbl{The implementation details  are discussed in the following \abrev{subsections}.}\footnote{\rbl{The 
module $\mathsf{SOLVE}$ is implemented for both $P_1$ and $P_2$ 
finite element \abrev{approximations} in the current version of the software,
whereas the error estimation module (and hence, the adaptive algorithm) 
is only implemented  for $P_1$ approximation.}}

\subsection{Module $\mathsf{SOLVE}$. Linear algebra aspects of  SGFEM} \label{subsec:asgfem:solve}	

Recalling the definitions of $\XX_h$ and $\PP_\gotP$,
the Galerkin solution $u_{\XX\PP}$ is sought in the form
\be \label{eq:approx:sol}
       u_{\XX\PP} = \sum_{i=1}^{N_\XX} \sum_{j=1}^{N_\PP} u_{ij} \phi_{i}(\bx) \sfP_{\kappa(j)}(\boldy),
\ee
where $N_\PP := \dim(\PP) = \# \gotP$, $\kappa$ is a bijection $\{1,2,\ldots,N_\PP\} \to \gotP$, and
the coefficients $u_{ij}$ are computed by solving the linear system $A \bu = \boldb$ with block structure.
Specifically, the solution vector $\bu$ and the right-hand side vector $\boldb$ are given by
\[
   \bu=[\bu_{1}\ \bu_{2}\ \ldots\ \bu_{N_{\PP}}]^{T} \text{ \ and \ }
   \boldb=[\boldb_{1}\ \boldb_{2}\ \ldots\ \boldb_{N_{\PP}}]^{T},
\]
respectively, with
\[
   \bu_{j} :=[u_{1j}\; u_{2j}\; \ldots\; u_{N_{\XX}j}]^{T},\quad j=1,\ldots,N_{\PP},
\]
\[
   [\boldb_{t}]_{s} := \langle 1,\sfP_{\kappa(t)} \rangle_\pi\, \int_{D} f(\bx) \phi_{s}(\bx)\, \dx,\quad
   s=1,\ldots,N_{\XX},\ \ t=1,\ldots,N_{\PP};
\]
the coefficient matrix $A$ is given by (see, e.g., \cite[Section~9.5]{lps2014})
\be \label{eq:decomp:A}
      A = G_{0} \otimes K_{0} + \sum_{m=1}^{M_\gotP} G_{m} \otimes K_{m},
\ee
where $M_\gotP$ is the number of active parameters in the index set $\gotP$,
\[
   [G_0]_{tj} := \langle \sfP_{\kappa(j)}, \sfP_{\kappa(t)} \rangle_{\pi} = \delta_{tj},\quad
   [G_m]_{tj} := \langle y_m \sfP_{\kappa(j)}, \sfP_{\kappa(t)} \rangle_{\pi}\ \ 
   (m=1,\ldots,M_\gotP)
\]
with $t,j = 1,\ldots,N_{\PP}$, and
$K_m$ are the finite element (stiffness) matrices defined by
\[
   [K_m]_{si} :=\int_{D} a_{m} \nabla\phi_{i} \cdot \nabla\phi_{s}\, \dx,\quad
   m=0,1,\ldots,M_\gotP,\ \
   s,i=1,\ldots,N_{\XX}.
\]

\rbl{The design of an efficient linear solver is a crucial ingredient of the stochastic
 Galerkin approximation process. Rather than computing a (memory intensive)
 sparse factorization of the coefficient matrix, a {\it matrix-free\/} iterative solver
 is needed.
 The key idea is that the \abrev{matrix-vector products with $A$ can be 
 computed  from its sparse matrix components  by
 exploiting the Kronecker product  structure, without assembling $A$ itself.}
 The  iterative solver that enables this process  within  T-IFISS is a 
 bespoke implementation of the Minimum Residual algorithm, called EST\_MINRES~\cite{ss11}. 
 The MINRES  algorithm is designed  to solve symmetric (possibly indefinite) linear equation systems
 and requires the action of $A$ on a given vector at every iteration, see~\cite[Section~2.4]{esw14}.
 Using  this strategy the storage overhead (in addition to the component matrices 
 $K_0, \ldots, K_{M_\gotP}, G_1, \ldots, G_{M_\gotP}$) is for
  five vectors of length $N_{\PP}\cdot N_{\XX}$.}
 
\rbl{A crucial ingredient in the design of a {\it fast} iterative solver is {\it preconditioning}. The
standard choice of  preconditioning operator  in this context is the 
parameter-free matrix operator 
$$P= G_{0} \otimes K_{0}=  I \otimes K_{0}. $$
The action of $P^{-1} \boldr$, where $\boldr$ is the current residual vector,
is needed at every iteration---this can be done efficiently
by computing  a single sparse triangular factorization of the matrix $K_0$ and then performing 
$N_\PP$  forward and backward substitutions on the components of the  residual vector.
Theoretical analysis of the preconditioned operator given in~\cite{pe2009} shows that the eigenvalues
of the preconditioned operator are bounded away from zero and bounded away from infinity  independently
of the \abrev{discretization} parameters  $N_{\XX}$ and  $N_{\PP}$. This means that the number of 
preconditioned EST\_MINRES iterations needed to satisfy a fixed residual
reduction tolerance will not grow unboundedly when the \abrev{discretization} parameters are changed. 
In practice, the number of iterations needed to satisfy the default tolerance of {\tt 1e-10} is less than 20, independent 
of the finite element mesh resolution  and the number of active indices.
}

\subsection{\rbl{Module $\mathsf{ESTIMATE}$. Error estimation in SGFEM}} \label{subsec:asgfem:estimate}

Stochastic Galerkin approximations are built from two distinct discretizations:
the spatial (finite element) discretization over the physical domain $D$ and the parametric (polynomial) approximation
on the parameter domain~$\G$.
Therefore, there are two distinct sources of discretization error arising from the choice
of the finite element space $\XX$ and the polynomial space $\PP$.
This fact determines the structure of a posteriori estimates for the energy errors in SGFEM approximations
as combinations of spatial and parametric contributions (cf.~\cite{egsz14, egsz15, bps14, bs16}).

The \emph{spatial} errors in SGFEM approximations are estimated by extending the strategies (EES1)--(EES3)
described in~\S\ref{sec:adaptivity} to tensor-product discretizations.
For example, in (EES1), each local (elementwise) error estimator, denoted by $e_{\XX}|_K$ ($K \in \mesh_h$),
now lives in the tensor-product space $\YY|_K \otimes \PP_\gotP$, where
$\YY|_K$ is the local space of piecewise linear or piecewise quadratic bubble functions.
These error estimators are computed by solving local residual problems of the following type (see~\cite[Section~6.2]{bps14} for details):
\be \label{eq:eX:local}
   B_{0,K}(e_{\XX}|_K,v) = \text{Res}_K(a,f,u_{\XX\PP};v)\quad\forall\,v \in \YY|_K \otimes \PP_\gotP,
\ee
where $B_{0,K}$ is the elementwise bilinear form associated with the parameter-free term $a_0$ in the coefficient expansion (cf.~\eqref{eq:B0}).
This construction of the error estimator enables fast linear algebra for solving~\eqref{eq:eX:local}.
Indeed, the coefficient matrix in the linear system associated with \eqref{eq:eX:local} has a very simple structure:
it is the Kronecker product of a $3 \times 3$ reduced stiffness matrix
and the identity matrix of dimension $N_\PP = \dim(\PP)$. 
As a result, the action of the inverse of this coefficient matrix
can be effected by a block $LDL^T$ factorization of the element stiffness matrices 
followed by a sequence of  $N_\PP$ backward and forward substitutions.
Furthermore, since the factorizations and triangular solves are logically independent,
the entire computation is vectorized over the finite elements that define the spatial subdivision.
We refer to~\cite[Section~3]{bprr19} for details of the global hierarchical (EES2)
and the two-level (EES3) error estimation strategies in the context of the SGFEM.

The \emph{parametric} errors in SGFEM approximations are estimated using the hierarchical approach in the spirit of~\cite{BankW_85_SAE}.
To that end, we first introduce the finite index set $\gotQ$
as a ``neighborhood'' of the index set $\gotP$.
More precisely, for a fixed $\overline{M} \in \NN$, we define
\begin{equation} \label{eq:index:set:Q}
   \gotQ := \big\{ \nu \in \gotI\setminus \gotP;\;
                              \nu  = \mu \pm \varepsilon^{(m)} \text{ for some } \mu \in \gotP
                              \text{ and some } m = 1,\dots, M_\gotP + \overline{M} \big\},
\end{equation}
where $\varepsilon^{(m)} := ( \varepsilon^{(m)}_1,  \varepsilon^{(m)}_2,\dots)$ ($m \in \NN$) denotes
the Kronecker delta index such that $\varepsilon^{(m)}_k =\delta_{mk}$ for all $k \in \NN$,
and $M_\gotP \in \NN$ is the number of active parameters in $\gotP$.

For a given $\gotP  \subset \gotI$, the index set $\gotQ$ contains
only those ``neighbors'' of all indices in $\gotP$ that have up to $M_\gotP + \overline{M}$ active parameters,
that is $\overline{M}$ parameters more than currently activated in the index set $\gotP$
(we refer to~\cite[Section~4.2]{bs16} for theoretical underpinnings of this construction).
Then, the \emph{parametric} error estimator $e_\PP$ is computed as a combination of the contributing estimators
$e_{\PP}^{(\nu)}$ associated with individual indices $\nu \in \gotQ$, i.e., $e_{\PP} = \sum_{\nu \in \gotQ} e_{\PP}^{(\nu)}$,
where each contributing estimator $e_{\PP}^{(\nu)} \in \XX \otimes \Span(\sfP_\nu)$, $\nu \in \gotQ$,
is computed by solving the linear system associated with the following discrete formulation:
\be \label{eq:e_Xnu}
      B_0(e_{\PP}^{(\nu)} , v \sfP_\nu) = F(v \sfP_\nu) - B(u_{\XX\PP}, v \sfP_\nu) \quad \text{for all $v \in \XX$}.
\ee
Note that the coefficient matrix of this linear system represents
the assembled stiffness matrix corresponding to the parameter-free term $a_0$ in~\eqref{eq1:a},
and is therefore the same for all $\nu \in \gotQ$.
Once the stiffness matrix has been factorized, the estimators $ e_\PP^{(\nu)}$
are computed independently by using forward and backward substitutions.

Once the \emph{spatial} and \emph{parametric} error estimators have been computed,
the total error estimate $\eta$ is calculated via
\be \label{eq:eta}
   \eta := \Big(\|e_{\XX}\|^2_{0} + \|e_{\PP}\|^2_{0}\Big)^{1/2} =
              \bigg( \sum\limits_{K \in \mesh_h} \big\|e_{\XX}|_{K}\|^2_{0,K} +
              \sum\limits_{\nu \in \gotQ} \big\|e_\PP^{(\nu)}\big\|^2_{0} \bigg)^{1/2},
\ee
where $\|\cdot\|_{0}$ (resp., $\|\cdot\|_{0,K}$) denotes the norm induced by the bilinear form $B_0$ (resp.,~$B_{0,K}$).

\subsection{Marking and refinement in adaptive SGFEM} \label{subsec:asgfem:adaptivity}

The module $\mathsf{ESTIMATE}$ supplies local spatial error indicators
associated with elements or edges of triangulation (e.g., $e_\XX|_K$ for the error estimation strategy (EES1)) as well as
the contributing parametric error indicators $e_{\PP}^{(\nu)}$ associated with individual indices $\nu \in \gotQ$.
In the module $\mathsf{MARK}$, the largest error indicators are selected independently for spatial 
and for parametric components of Galerkin approximations.
To that end, one of the marking strategies described in~\S\ref{sec:adaptivity}
(i.e., either the maximum or the D\"orfler strategy) is employed.
In Stochastic T-IFISS, the same marking strategy is used for both spatial and parametric components
with marking thresholds $\theta_\XX$ and $\theta_\PP$, respectively.
However, a simple modification of the code will allow one to use different marking strategies for different 
components of Galerkin approximations.

Thus, at each iteration of the adaptive SGFEM algorithm, the output of the module $\mathsf{MARK}$ contains two sets:
the set of marked elements in the current mesh $\mesh_h$ to be refined (or, the set of edges to be bisected)
and the set $\gotM \subseteq \gotQ$ of marked indices to be added to the current index set $\gotP$
(note that choosing $\overline{M} > 1$ in~\eqref{eq:index:set:Q} allows one to activate more than one
new parameter at the next iteration of the adaptive loop).
The finite-dimensional space $V_{\XX\PP}$ is then enhanced within the module $\mathsf{REFINE}$ by performing
either spatial refinement (as described in~~\S\ref{sec:adaptivity}) or
parametric refinement (simply by adding $\gotM$ to $\gotP$).
The question then arises which type of refinement (spatial {\it or} parametric) should be performed at a given iteration.

A traditional strategy for choosing between the two refinements is based on the dominant error estimator contributing 
to the {total} error estimate $\eta$ defined in~\eqref{eq:eta}; cf.~\cite{egsz14, egsz15, bs16}.
This strategy is referred to as \emph{version~1} of the adaptive algorithm implemented in Stochastic T-IFISS.
An alternative strategy is referred to as \emph{version~2} of the implemented algorithm:
here, the refinement type that leads to a larger estimated error reduction is chosen at each
 iteration; see~\cite{br18, BespalovPRR_CAS}.
This strategy exploits the fact that local \emph{spatial} error indicators
(e.g., $\|e_\XX|_K\|_{0,K}$ ($K \in \mesh_h$) in the error estimation strategy (EES1)) and
individual \emph{parametric} error indicators $\|e_{\PP}^{(\nu)}\|_0$ ($\nu \in \gotQ$)
provide effective estimates of the error reduction that would be achieved by performing, respectively,
a local refinement of the current mesh (e.g., by refining the element $K$) and
a selective enrichment of the parametric component of the current Galerkin approximation
(by adding the index $\nu \in \gotQ$ to the current index set $\gotP$).
We refer to~\cite[Theorem~5.1]{bps14} and~\cite[Corollary~3]{BespalovPRR_CAS} for the 
underpinning theoretical results
and to~\cite{br18} and~\cite[Section~5]{BespalovPRR_CAS} for comprehensive
numerical studies of the two versions of the adaptive algorithm and different marking strategies.


\subsection{\rev{Numerical case study}} \label{subsec:asgfem:numerics}


\rev{We conclude this section  with a representative case study that demonstrates the efficiency of our
 adaptive SGFEM algorithm.}

\smallskip

\begin{figure}[t!] 
\centering
\subfloat[][]
{
\includegraphics[scale=0.93]{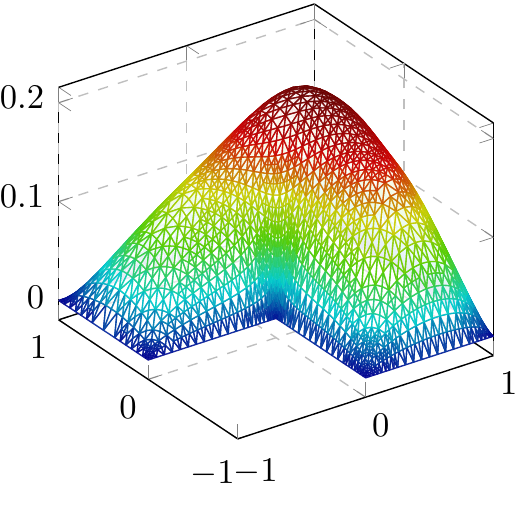}}   
\hfill
\subfloat[][] 
{
\includegraphics[scale=0.93]{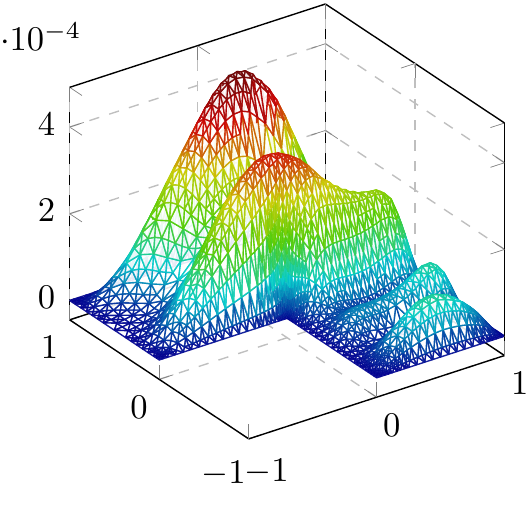}}  
\hfill 
\subfloat[][]
{
\includegraphics[scale=0.93]{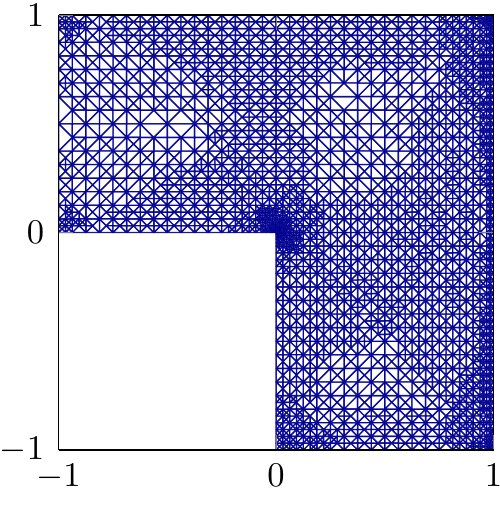}}  
\caption{
\rbl{Example 4}: 
(a)--(b) mean and variance of the SGFEM solution;
(c) a typical locally refined mesh generated by the adaptive SGFEM algorithm. 
}
\label{fig:ex5:mean:var:mesh}
\end{figure}

{\bf Example~4.}
We consider the parametric model problem~\eqref{eq:strongform} on the L-shaped domain 
$D = (-1,1)^2\setminus (-1,0]^2$.
We set $f(\bx) = (1-x_1)^{-0.4}$
and choose the parametric coefficient $a(\bx,\boldy)$ in the form~\eqref{eq1:a},
where $a_0$ and $a_m$ ($m \in \NN$) are as specified by the coefficient expansion (CE2) 
with $\tilde\sigma = 2$
and $y_m$ ($m\in \NN$) are the images of
independent truncated Gaussian random variables with zero mean on 
$\Gamma_m = [-1,1]$ (see (RV2), where we set $\sigma_0 = 1$).
The mean and the variance of an SGFEM solution to this problem are shown in 
Figure~\ref{fig:ex5:mean:var:mesh}(a)--(b),
whereas Figure~\ref{fig:ex5:mean:var:mesh}(c) depicts a typical locally refined mesh 
generated by the adaptive SGFEM algorithm.
Note how the adaptively refined mesh effectively identifies the areas of singular
behavior of the mean field---in the
vicinity of the reentrant corner (due to a geometric singularity) and in the vicinity of the edge 
$x_1 = 1$ (due to a singular right-hand side function).

When running the adaptive SGFEM algorithm, 
\rbl{we start with a uniform mesh consisting of \abrev{96} right-angled triangles
and an initial index set} $\gotP_0 := \{(0,0,0,\ldots,),\; (1,0,0,\ldots,)\}$.
For the error estimation module, we employ the two-level spatial error estimator (EES3)
combined with the hierarchical parametric error estimators
associated with individual indices $\nu \in \gotQ$ (see~\eqref{eq:e_Xnu} and~\eqref{eq:index:set:Q},
 where we set $\overline{M} := 1$).
We use D\"orfler marking for spatial and parametric components of Galerkin approximations
(for the spatial component, we mark the edges of the mesh).

\begin{figure}[b!]
\centering
\begin{tikzpicture}
\pgfplotstableread{data/exp5/thetaX1.0_thetaP1.0.dat}{\oneone}
\pgfplotstableread{data/exp5/thetaX0.7_thetaP1.0.dat}{\zerosevenone}
\pgfplotstableread{data/exp5/thetaX1.0_thetaP0.7.dat}{\onezeroseven}
\pgfplotstableread{data/exp5/thetaX0.7_thetaP0.9.dat}{\zerosevenzeronine}
\begin{loglogaxis}
[
width = 10cm, height=7cm,									
xlabel={degrees of freedom},						
ylabel={error estimate},							
ymajorgrids=true, xmajorgrids=true, grid style=dashed,		
xmin = (8)*10^(0), 	xmax = (8)*10^(7), 
ymin = (5)*10^(-4), ymax = (3)*10^(-1),					
legend style={legend pos=south west, legend cell align=left, fill=none, draw=none, font={\fontsize{9pt}{12pt}\selectfont}}
]
\addplot[red,mark=o,mark size=2.5pt]			table[x=dofs, y=error]{\oneone};
\addplot[myOrange,mark=pentagon,mark size=3.0pt]	table[x=dofs, y=error]{\zerosevenone};
\addplot[darkGreen,mark=square,mark size=2.0pt]			table[x=dofs, y=error]{\onezeroseven};
\addplot[blue,mark=triangle,mark size=3.5pt]		table[x=dofs, y=error]{\zerosevenzeronine};
\addplot[black,solid,domain=8*10^(0):8*10^(7)]    { (0.4)*x^(-1/3) };
%
\legend{
{$\theta_\XX=1$, $\theta_\PP=1$},
{$\theta_\XX=0.7$, $\theta_\PP=1$},
{$\theta_\XX=1$, $\theta_\PP=0.7$},
{$\theta_\XX=0.7$, $\theta_\PP=0.9$},
{$\mathcal{O}(N^{-1/3})$}
}
\end{loglogaxis}
\end{tikzpicture}
\caption{
\rbl{Example 4}: error estimates at each iteration of the adaptive SGFEM algorithm
for different sets of marking parameters.
}
\label{fig:ex5:adaptivity:wins}
\end{figure}
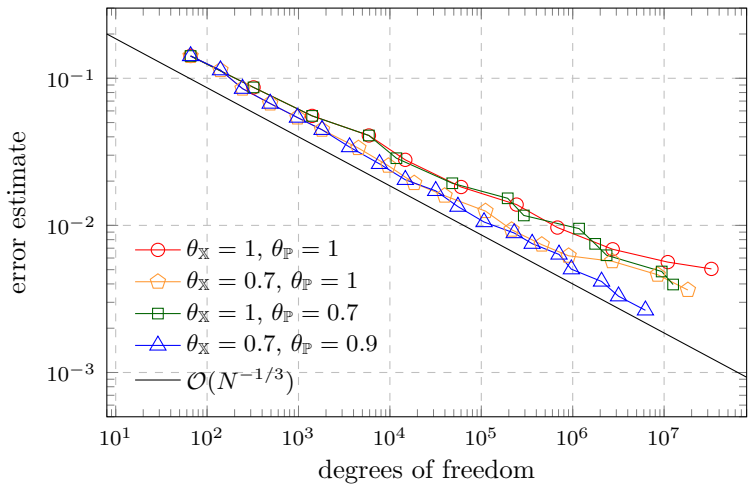

\rbl{In our first experiment we ran} the adaptive SGFEM algorithm (in version 2) 
with the following four sets of D\"orfler marking parameters
(\rbl{with} different stopping tolerances):\\
(i) $\theta_\XX = 1$, $\theta_\PP = 1$ (no adaptivity in either of the components),  {\tt tol=5.1e-3};\\
(ii) $\theta_\XX = 0.7$, $\theta_\PP = 1$ (adaptive refinement only for spatial component),  {\tt tol=4e-3};\\
(iii) $\theta_\XX = 1$, $\theta_\PP = 0.7$ (adaptive refinement only for parametric component), {\tt tol=4e-3} ;\\
(iv) $\theta_\XX = 0.7$, $\theta_\PP = 0.9$ (adaptive refinement of both components),  {\tt tol=3e-3}.\\
For each run of the algorithm, the error estimates computed at each iteration are plotted in 
Figure~\ref{fig:ex5:adaptivity:wins}.
\rbl{The error estimates can be seen to decrease at every iteration.}
However, in cases (i)--(iii), the decay rate either eventually 
deteriorates (cases (i) and (ii)), due to the number of degrees of freedom
growing very fast, or it is significantly slower (case (iii)) than in the case of adaptivity being used
for \emph{both} components of SGFEM approximations (case (iv)).
This shows that, for the same level of accuracy, adaptive refinement of \emph{both} components
results in more balanced approximations with \rev{fewer} degrees of freedom and leads to 
\rbl{the fastest convergence rate for the parameter choices  in this experiment}.

\begin{table}[!t]
\begin{center}
\begin{tabular}{ccc}
 Iteration &  & Evolution of the index set                    \\[3pt]
  0	&	& (0 0 0 0 0 0)      \\
	&	& (1 0 0 0 0 0)      \\[3pt]
  7	&	& (0 1 0 0 0 0)      \\
	&	& (2 0 0 0 0 0)      \\[3pt]
  11	&	& (0 0 1 0 0 0)      \\
	&	& (1 1 0 0 0 0)      \\
	&	& (3 0 0 0 0 0)      \\[3pt]
  14	&	& (0 0 0 1 0 0)      \\
	&	& (1 0 1 0 0 0)      \\
	&	& (2 1 0 0 0 0)      \\
	&	& (0 2 0 0 0 0)      \\[3pt]
  16	&	& (0 0 0 0 1 0)      \\
	&	& (2 0 1 0 0 0)      \\
	&	& (1 0 0 1 0 0)      \\
	&	& (4 0 0 0 0 0)      \\[3pt]
  18	&	& (0 0 0 0 0 1)      \\
	&	& (1 0 0 0 1 0)      \\
	&	& (3 1 0 0 0 0)      \\
	&	& (0 1 1 0 0 0)      \\
	&	& (1 2 0 0 0 0)      \\
	&	& (2 0 0 1 0 0)      \\
	&	& (3 0 1 0 0 0)      \\
	&	& (1 0 0 0 0 1)
\end{tabular}
\caption{\rbl{Example 4}: evolution of the index set when running adaptive SGFEM algorithm with
$\theta_\XX = 0.7$ \rbl{and $\theta_\PP = 0.9$}.
}
\label{tbl:ex5:indset}
\end{center}
\end{table}

\begin{figure}[t!]
\subfloat[][] 
{
\begin{tikzpicture}
\pgfplotstableread{data/exp5/thetaX0.7_thetaP0.9.dat}{\first}
\begin{loglogaxis}
[
title={Adaptive SGFEM ($\theta_\XX = 0.7$,  $\theta_\PP = 0.9$)},				
xlabel={degrees of freedom},							
ylabel={error estimate},								
ymajorgrids=true, xmajorgrids=true, grid style=dashed,			
xmin = 8*10^(0), 	xmax = 2*10^(7),								
ymin = 1*10^(-3),	ymax = 2*10^(-1),								
width = 7.5cm, height=7.0cm,
legend style={legend pos=south west, legend cell align=left, fill=none, draw=none, font={\fontsize{9pt}{12pt}\selectfont}}
]
\addplot[blue,mark=triangle,mark size=3pt]	table[x=dofs, y=error]{\first};
\addplot[darkGreen,mark=o,mark size=2pt]	table[x=dofs, y=yp_one]{\first};
\addplot[red,mark=square,mark size=2pt]	table[x=dofs, y=xq_one]{\first};
\legend{
{$\eta$ (total)},
{$\|e_\XX\|_{0}$ (spatial)},
{$\|e_\PP\|_{0}$ (parametric)},
}
\end{loglogaxis}
\end{tikzpicture}
}
\hfill
\subfloat[][] 
{
\begin{tikzpicture}
\pgfplotstableread{data/exp5/thetaX0.7_thetaP0.9.dat}{\first}
\begin{semilogxaxis}
[
title={Adaptive SGFEM ($\theta_\XX = 0.7$,  $\theta_\PP = 0.9$)},				
xlabel={degrees of freedom}, 						
ylabel={effectivity index},								
ymajorgrids=true, xmajorgrids=true, grid style=dashed,		
xmin = 8*10^(0), 	xmax = 2*10^(7),								
ymin = 0.6,	ymax = 1.0,						
width = 7.5cm, height=7.0cm,
ytick={0.6,0.7,0.8,0.9,1.0},
]
\addplot[blue,mark=triangle,mark size=3.0pt]	table[x=dofs, y=effindices]{\first};
\end{semilogxaxis}
\end{tikzpicture}
}
\caption{
\rbl{Example 4}: (a) energy error estimates at each iteration, along with their spatial and parametric contributions;
(b) the associated effectivity indices.
}
\label{fig:ex5:estimates:effectivities}
\end{figure}
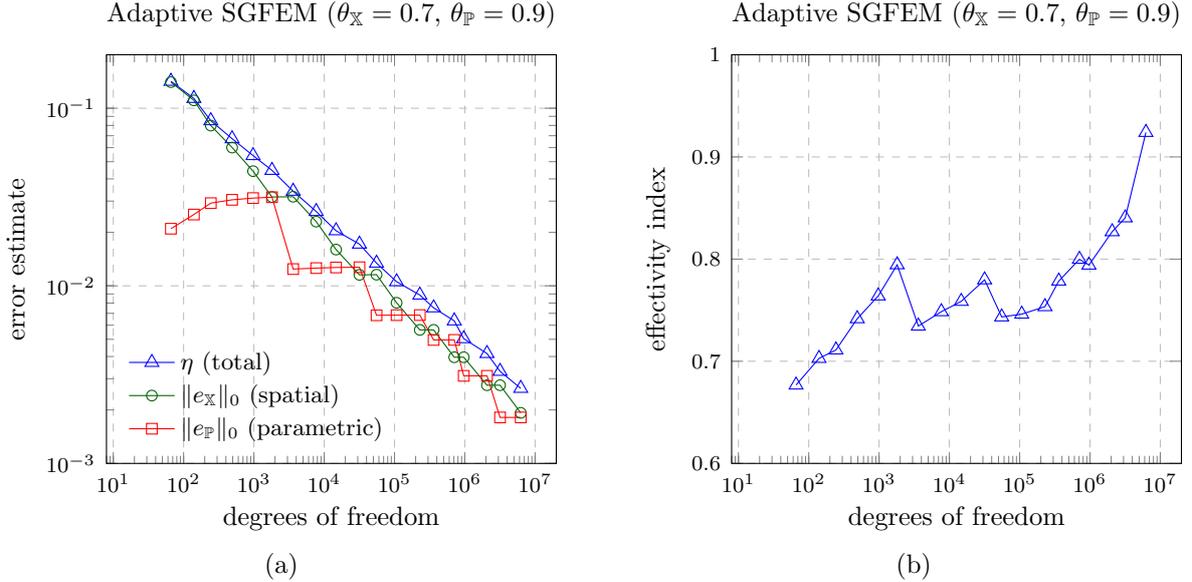

\rbl{To give an indication of the algorithmic efficiency, we provide further details of the run} in case~(iv).
\rbl{In this computation}, which took 927 seconds, the stopping tolerance  {\tt 3e-3} was met after 19 iterations.
The final triangulation generated by the algorithm comprised 545,636 finite elements with 271,599 interior vertices
(the latter number defines the dimension of the corresponding finite element space).
For the final polynomial approximation on $\G$, the algorithm produced an index set $\gotP$ of 
cardinality 23 with 6 active parameters;
the evolution of the index set throughout the computation is shown in Table~\ref{tbl:ex5:indset}.
The total number of d.o.f. in the SGFEM solution at the final iteration was equal to 6,246,777.
In Figure~\ref{fig:ex5:estimates:effectivities}(a), we show the interplay between the spatial and parametric contributions
to the total error estimates $\eta$ (see~\eqref{eq:eta}) at each iteration of the adaptive algorithm.
The effectivity indices for the total error estimates are plotted in Figure~\ref{fig:ex5:estimates:effectivities}(b).
These were calculated using the reference Galerkin solution
computed with $P_2$ approximations on the final triangulation generated by the algorithm
and with the polynomial space $\PP_{\gotP \cup \gotQ}$, where $\gotP$ is the final index set produced by the algorithm
and $\gotQ$ is the ``neighborhood'' of $\gotP$ as defined in~\eqref{eq:index:set:Q} with $\overline{M} = 1$
(the total number of d.o.f. in this reference solution was~37,020,322).

